\documentclass[epsfig,latexsym,amsfonts,twoside]{article}
\usepackage{amssymb,amsmath,amscd}

\def\part#1{\frac{\partial\phantom{#1}}{\partial#1}}
\newtheorem{thm}{Theorem}
\newtheorem{prp}[thm]{Proposition}
\newtheorem{lem}[thm]{Lemma}

\newtheorem{cnj}[thm]{Conjecture}
\newenvironment{prf}{\begin{trivlist}\item[]{\bf Proof} }%
{\hfill $\Box$ \end{trivlist}}
\newenvironment{dfn}{\begin{trivlist}\item[]{\bf Definition}\em }%
{\end{trivlist}}
\newenvironment{rmk}{\begin{trivlist}\item[]{\bf Remark} }%
{\end{trivlist}}
\newenvironment{exm}{\begin{trivlist}\item[]{\bf Example} }%
{\end{trivlist}}

\def\Z{\ifmmode{{\mathbb Z}}\else{${\mathbb Z}$}\fi}
\def\Q{\ifmmode{{\mathbb Q}}\else{${\mathbb Q}$}\fi}
\def\C{\ifmmode{{\mathbb C}}\else{${\mathbb C}$}\fi}
\def\P{\ifmmode{{\mathbb P}}\else{${\mathbb P}$}\fi}

\def\H{\ifmmode{{\mathrm H}}\else{${\mathrm H}$}\fi}

\def\B{\ifmmode{{\cal B}}\else{${\cal B}$}\fi}
\def\E{\ifmmode{{\cal E}}\else{${\cal E}$}\fi}
\def\F{\ifmmode{{\cal F}}\else{${\cal F}$}\fi}
\def\I{\ifmmode{{\cal I}}\else{${\cal I}$}\fi}
\def\K{\ifmmode{{\cal K}}\else{${\cal K}$}\fi}
\def\L{\ifmmode{{\cal L}}\else{${\cal L}$}\fi}
\def\M{\ifmmode{{\cal M}}\else{${\cal M}$}\fi}
\def\N{\ifmmode{{\cal N}}\else{${\cal N}$}\fi}
\def\O{\ifmmode{{\cal O}}\else{${\cal O}$}\fi}
\def\U{\ifmmode{{\cal U}}\else{${\cal U}$}\fi}
\def\X{\ifmmode{{\cal X}}\else{${\cal X}$}\fi}

\def\Br{\ifmmode{{\mathrm{Br}}}\else{${\mathrm{Br}}$}\fi}
\def\OG{\ifmmode{\widetilde{\cal M}_4}\else{$\widetilde{\cal M}_4$}\fi}
\def\D{\ifmmode{{\cal D}_{\mathrm{coh}}^b}\else{${{\cal
    D}_{\mathrm{coh}}^b}$}\fi}
\def\Shah{\ifmmode{\amalg\hspace*{-3.5pt}\amalg}\else{$\amalg\hspace*{-3.5pt}\amalg$}\fi}

\begin{document}

\title{On Lagrangian fibrations by Jacobians II}
\author{Justin Sawon}
\date{September, 2011}
\maketitle

\begin{abstract}
Let $Y\rightarrow\P^n$ be a flat family of reduced Gorenstein curves, such that the compactified relative Jacobian $X=\overline{J}^d(Y/\P^n)$ is a Lagrangian fibration. We prove that $X$ is a Beauville-Mukai integrable system if $n=3$, $4$, or $5$, and the curves are irreducible and non-hyperelliptic. We also prove that $X$ is a Beauville-Mukai system if $n=3$, $d$ is odd, and the curves are canonically positive $2$-connected hyperelliptic curves.
\end{abstract}

\section{Introduction}

In this article a {\em Lagrangian fibration\/} will be a holomorphic symplectic manifold $X$ of dimension $2n$, with a proper surjective morphism $\pi:X\rightarrow\P^n$. The fibres of $\pi$ will be Lagrangian with respect to the holomorphic symplectic form on $X$, and the generic fibre will be an $n$-dimensional abelian variety (see Matsushita~\cite{matsushita99, matsushita00i}). The following conjecture was introduced in Part~I~\cite{sawon08iv} of this article.
\begin{cnj}
\label{bmsystem}
Let $X\rightarrow\P^n$ be a Lagrangian fibration whose fibres are Jacobians of genus $n$ curves. Then $X$ is a Beauville-Mukai integrable system.
\end{cnj}
A description of the Beauville-Mukai systems~\cite{beauville99, mukai84} was given in Part~I. Essentially it means that $X$ is the compactified relative Jacobian of a complete linear system of curves in a K3 surfaces.

To be more precise, suppose that $Y\rightarrow\P^n$ is a flat family of integral Gorenstein curves of arithmetic genus $n$, whose compactified relative Jacobian
$$X=\overline{J}^d(Y/\P^n)$$
is a Lagrangian fibration. When $n=2$, Markushevich~\cite{markushevich96} proved that $X$ must be a Beauville-Mukai system. In Part~I we extended this result to arbitrary dimension, provided the degree of the discriminant locus $\Delta\subset\P^n$ is greater than $4n+20$.

In this article, we prove Conjecture~\ref{bmsystem} in various low dimensional cases. Firstly, we prove $X$ is a Beauville-Mukai system if all of the curves in the family $Y/\P^n$ are non-hyperelliptic and $n=3$, $n=4$, or $n=5$ (Theorems~\ref{g=3nh},~\ref{g=4}, and~\ref{g=5}). Secondly, when $n=3$ we prove that the family of curves $Y/\P^3$ cannot contain both non-hyperelliptic curves and hyperelliptic curves (Proposition~\ref{notboth}). This is in accordance with the well-known fact that given a linear system of curves in a K3 surface, either all of the curves are hyperelliptic or none of the curves are hyperelliptic. Thirdly, when $n=3$ we prove that $X$ is a Beauville-Mukai system if all of the curves in the family $Y/\P^3$ are hyperelliptic (Theorem~\ref{g=3h}). Note that the hypotheses need to be modified slightly in this case, as we cannot expect all of the curves in the family to be irreducible. Instead we require the curves to be canonically positive and $2$-connected, which ensures that the curves do not contain `redundant' irreducible components. In addition, the compactified relative Jacobians of different degrees may not be locally isomorphic as fibrations if there are reducible curves; our result assumes that the degree $d$ is odd.

Unfortunately our methods do not immediately extend to higher genus, since we use the fact that the curves are complete intersections. In the non-hyperelliptic case, the canonical embedding of a generic curve of genus $g\geq 6$ does not yield a complete intersection in $\P^{g-1}$. In the hyperelliptic case, the image of the canonical map of a curve of genus $g\geq 4$ yields a rational normal curve in $\P^{g-1}$, which again is not a complete intersection.

The author would like to thank Eduardo Esteves and Steven Kleiman for patiently answering many of his questions concerning compactified Jacobians, and would also like to thank the Max-Planck-Institut f{\"u}r Mathematik in Bonn, where these results were obtained.

\section{Preliminaries}

\subsection{Canonical models of curves}

The following results were proved in Sections~2.2 and 4.2 of Part I~\cite{sawon08iv}.

\begin{lem}
\label{preliminaries}
Let $Y\rightarrow B$ be a flat family of geometrically integral (reduced and irreducible) curves over a projective manifold $B$, such that the compactified relative Jacobian
$$X=\overline{J}^d(Y/B)$$
is a Lagrangian fibration over $B$. Then
\begin{enumerate}
\item every curve in the family $Y/B$ has arithmetic genus $n$,
\item the base $B$ is isomorphic to $\P^n$,
\item the generic curve in the family $Y/B$ is a smooth genus $n$ curve,
\item there is a hypersurface $\Delta\subset B$ parametrizing singular fibres of $\pi:X\rightarrow B$ (equivalently, singular curves in the family $Y/B$) and a curve above a generic point of $\Delta$ will contain a single simple node,
\item the compactified relative Jacobian $X=\overline{J}^d(Y/B)$ is a Lagrangian fibration if and only if $X^1=\overline{J}^1(Y/B)$ is a Lagrangian fibration,
\item the first direct image sheaf $R^1\pi_*\O_X$ is isomorphic to $\Omega^1_{\P^n}$, where $\pi$ denotes the projection $X\rightarrow B\cong\P^n$,
\item the first direct image sheaf $R^1\pi_*\O_Y$ is isomorphic to $\Omega^1_{\P^n}$, where $\pi$ denotes the projection $Y\rightarrow B\cong\P^n$.
\end{enumerate}
If in addition the curves in the family $Y/B$ are Gorenstein, then the total space $Y$ of the family of curves is smooth.
\end{lem}

Suppose that we are in the above situation, i.e., suppose that $Y/\P^n$ is a flat family of integral Gorenstein curves of genus $n$, such that $X=\overline{J}^d(Y/\P^n)$ is a Lagrangian fibration.

The canonical map
$$C\rightarrow\P(\H^0(C,K_C)^{\vee})\cong\P(\H^1(C,\O_C))\cong\P^{n-1}$$
of a smooth curve of genus $n\geq 2$ is either an embedding or two-to-one onto a rational normal curve, according to whether $C$ is non-hyperelliptic or hyperelliptic (see Griffiths and Harris~\cite{gh78}). Rosenlicht~\cite{rosenlicht52} extended the canonical map to singular curves: he used $\H^0(\omega)$ to define a base-point-free linear series on the normalization $\tilde{C}$ of $C$, and defined the canonical model of $C$ to be the image $C^{\prime}$ of the corresponding map $\tilde{C}\rightarrow\P^{n-1}$. He then proved the following result (see Kleiman and Martins~\cite{km09} for a modern presentation).

\begin{thm}
\label{rosenlicht}
Let $C^{\prime}$ be the canonical model of an integral curve $C$ of genus $n$.
\begin{enumerate}
\item If $C$ is not birational to $C^{\prime}$, then $C^{\prime}$ is a rational normal curve of degree $n-1$ in $\P^{n-1}$ and the map $\tilde{C}\rightarrow\P^{n-1}$ induces a unique two-to-one map $C\rightarrow\P^1$. We call such a curve hyperelliptic. Note that $C$ is automatically Gorenstein in this case.
\item If $C$ is non-hyperelliptic and Gorenstein, then $C^{\prime}$ is degree $2n-2$ in $\P^{n-1}$ and the map $\tilde{C}\rightarrow\P^{n-1}$ induces an isomorphism $C\cong C^{\prime}$.
\item If $C$ is non-hyperelliptic and non-Gorenstein, then $C^{\prime}$ is isomorphic to the blowup of $C$ with respect to $\omega$.
\end{enumerate}
\end{thm}

By hypothesis, the curves in our family $Y/\P^n$ are Gorenstein. So case (3) of the above theorem does not arise, and there exists a relative canonical map
$$Y\rightarrow\P(R^1\pi_*\O_Y)\cong\P(\Omega^1_{\P^n})$$
where we have used Lemma~\ref{preliminaries}~(7). We will treat the non-hyperelliptic and hyperelliptic cases separately. But first we study the geometry of $\P(\Omega^1_{\P^n})$.

\subsection{The incidence variety}

\begin{dfn}
The incidence variety is the subvariety $V\subset\P^n\times(\P^n)^{\vee}$ parametrizing pairs $(p,H)$, where $p$ is a point in the hyperplane $H$.
\end{dfn}

It is easily verified that $V$ is a hypersurface of bidegree $(1,1)$. It is a $\P^{n-1}$-bundle over both $\P^n$ and $(\P^n)^{\vee}$ under the projections
$$\begin{array}{ccccc}
 & & V\subset{\P}^n\times({\P}^n)^{\vee} & & \\
 & h\swarrow & & \searrow j & \\
 {\P}^n & & & & ({\P}^n)^{\vee}.
 \end{array}$$
We prove some basic facts about $V$.

\begin{lem}
\label{aboutV}
\begin{enumerate}
\item Restriction of divisors induces an isomorphism
$$\mathrm{Pic}V\cong\mathrm{Pic}(\P^n\times(\P^n)^{\vee})\cong\Z\oplus\Z.$$
\item For $n\geq2$, the line bundle $\O(a,b)|_V=h^*\O_{\P^n}(a)\otimes j^*\O_{(\P^n)^{\vee}}(b)$ admits non-trivial sections if and only if $a\geq 0$ and $b\geq 0$. Every section lifts to a section of $\O(a,b)$ on $\P^n\times(\P^n)^{\vee}$.
\item For $n\geq 3$, $\H^1(V,\O(a,b)|_V)$ is trivial for all $a$ and $b$.
\end{enumerate}
\end{lem}

\begin{prf}
The first statement follows from the Lefschetz Hyperplane Theorem. Next consider the short exact sequence
$$0\rightarrow\O(a-1,b-1)\rightarrow\O(a,b)\rightarrow\O(a,b)|_V\rightarrow 0$$
and its corresponding long exact sequence
$$0\rightarrow\H^0(\O(a-1,b-1))\rightarrow\H^0(\O(a,b))\rightarrow\H^0(V,\O(a,b)|_V)\rightarrow\H^1(\O(a-1,b-1))$$
$$\rightarrow\H^1(\O(a,b))\rightarrow\H^1(V,\O(a,b)|_V)\rightarrow\H^2(\O(a-1,b-1))\rightarrow\ldots$$
A line bundle on $\P^n$ can only have non-vanishing cohomology in degree $0$ or $n$. Therefore on $\P^n\times(\P^n)^{\vee}$
$$\H^k(\O(c,d))=\bigoplus_{i+j=k}\H^i(\P^n,\O(c))\otimes\H^j((\P^n)^{\vee},\O(d))$$
is always trivial for $1\leq k\leq n-1$; this proves the second and third statements.
\end{prf}

We introduced the incidence variety because it is related to the $\P^{n-1}$-bundle $\P(\Omega^1_{\P^n})$.

\begin{lem}
\label{V}
The projective bundle $\P(\Omega^1_{\P^n}(1))$ is canonically isomorphic to $V$. Under this isomorphism, the relative hyperplane bundle $\O_{\P(\Omega^1_{\P^n}(1))}(1)$ is identified with $\O(0,1)|_V$.
\end{lem}

\begin{prf}
Projectivizing the Euler sequence
$$0\rightarrow\Omega^1_{\P^n}(1)\rightarrow\H^0(\P^n,\O(1))\otimes\O\stackrel{\mathrm{ev}}{\longrightarrow}\O(1)\rightarrow 0$$
we see that $\P(\Omega^1_{\P^n}(1))$ sits inside $\P^n\times(\P^n)^{\vee}$. Moreover, the inclusion of this $\P^{n-1}$-bundle into the trivial $(\P^n)^{\vee}$-bundle is compatible with the relative polarizations. This means that $\O_{\P(\Omega^1_{\P^n}(1))}(1)$ is the restriction of the relative hyperplane bundle $\O(0,1)$ on $\P^n\times(\P^n)^{\vee}$.

Finally, we identify $\P(\Omega^1_{\P^n}(1))$ with $V$. Suppose that the non-trivial section $s\in\H^0(\P^n,\O(1))$ vanishes on the hyperplane $H$. Then $(p,H)$ lies in $\P(\Omega^1_{\P^n}(1))$ if and only if $s(p)=0$, or equivalently, $p\in H$ and $(p,H)$ lies in $V$.
\end{prf}

\begin{rmk}
Lemma~\ref{V} implies there is also an isomorphism
$$\P(\Omega^1_{\P^n})\cong\P(\Omega^1_{\P^n}(1))=V.$$
Under this isomorphism, the relative hyperplane bundle $\O_{\P(\Omega^1_{\P^n})}(1)$ is identified with
$$\O(0,1)|_V\otimes h^*\O_{\P^n}(1)=\O(1,1)|_V.$$
\end{rmk}



\section{Non-hyperelliptic curves}

\subsection{Genus three}

\begin{thm}
\label{g=3nh}
Let $Y\rightarrow\P^3$ be a flat family of integral Gorenstein non-hyperelliptic curves of genus three whose compactified relative Jacobian $X=\overline{J}^d(Y/\P^3)$ is a Lagrangian fibration. Then $X$ is a Beauville-Mukai integrable system, i.e., the family of curves is a complete linear system of curves in a K3 surface.
\end{thm}

\begin{prf}
By Theorem~\ref{rosenlicht}~(2), the relative canonical map
$$Y\rightarrow \P(\Omega^1_{\P^3})\cong V\subset\P^3\times (\P^3)^{\vee}$$
is a closed embedding, with each curve $Y_t$ embedded as a quartic in the corresponding $\P^2$ fibre of $\P(\Omega^1_{\P^3})$. Regarding $Y$ as a subvariety of $V$, we see that $Y$ will be a hypersurface given by the zero locus of a section in
$$\H^0(\P(\Omega^1_{\P^3}),\O_{\P(\Omega^1_{\P^3})}(4)\otimes h^*\O_{\P^3}(d))=\H^0(V,\O(d+4,4)|_V)$$
for some integer $d$. By Lemma~\ref{aboutV}~(2), $d$ must be at least $-4$ for this space of sections to be non-trivial.

\vspace*{3mm}
\noindent
{\bf Claim:} $d=-4$.

\vspace*{3mm}
By Lemma~\ref{preliminaries}~(7) we know that $R^1\pi_*\O_Y\cong\Omega^1_{\P^3}$. Thinking of $Y$ as a subvariety of $V$, we compute this a second time using the projection $h:V\rightarrow\P^3$. We start with a resolution of $\O_Y$ by locally free sheaves:
$$0\rightarrow\O_{\P(\Omega^1_{\P^3})}(-4)\otimes h^*\O_{\P^3}(-d)\rightarrow\O_{\P(\Omega^1_{\P^3})}\rightarrow\O_Y\rightarrow 0$$
Since $R^1h_*\O_{\P(\Omega^1_{\P^3})}$ and $R^2h_*\O_{\P(\Omega^1_{\P^3})}$ both vanish, the long exact sequence obtained by applying $R^{\bullet}h_*$ to the above yields
\begin{eqnarray*}
R^1h_*{\O}_Y & \cong & R^2h_*\left({\O}_{{\P}(\Omega^1_{{\P}^3})}(-4)\otimes h^*{\O}_{{\P}^3}(-d)\right) \\
 & \cong & {\O}_{{\P}^3}(-d)\otimes R^2h_*\left({\O}_{{\P}(\Omega^1_{{\P}^3})}(-4)\right) \\
 & \cong & {\O}_{{\P}^3}(-d)\otimes \left(h_*\left({\O}_{{\P}(\Omega^1_{{\P}^3})}(4)\otimes\omega_h\right)\right)^{\vee} \\
 & \cong & {\O}_{{\P}^3}(-d)\otimes \left(h_*\left({\O}_{{\P}(\Omega^1_{{\P}^3})}(4)\otimes{\O}_{{\P}(\Omega^1_{{\P}^3})}(-3)\otimes h^*\omega_{{\P}^3}^{\vee}\right)\right)^{\vee} \\
 & \cong & {\O}_{{\P}^3}(-d)\otimes\omega_{{\P}^3}\otimes \left(h_*{\O}_{{\P}(\Omega^1_{{\P}^3})}(1)\right)^{\vee} \\
 & \cong & {\O}_{{\P}^3}(-d-4)\otimes \Omega^1_{{\P}^3}
 \end{eqnarray*}
where we have used relative Serre duality on the third line and the fact that
$$\omega_V\cong\omega_{\P^3\times(\P^3)^{\vee}}\otimes\O(1,1)|_V\cong\O(-3,-3)|_V\cong\O_{\P(\Omega^1_{\P^3})}(-3)$$
on the fourth line. Since $R^1h_*\O_Y\cong R^1\pi_*\O_Y\cong\Omega^1_{\P^3}$, $d$ must equal $-4$, proving the claim.

\vspace*{3mm}
We have proved that $Y\subset V$ is the zero locus of a section of $\O(0,4)|_V$. By Lemma~\ref{aboutV}~(2), this section lifts to a section of $\O(0,4)$ on $\P^3 \times (\P^3)^{\vee}$, which defines a quartic K3 surface $S$ in $(\P^3)^{\vee}$. The projection $j:V\rightarrow (\P^3)^{\vee}$ expresses $Y$ as a $\P^2$-bundle over $S$. In particular, $S$ must be smooth since $Y$ is smooth by the last part of Lemma~\ref{preliminaries}. 
$$\begin{array}{ccccccc}
 & & V \supset Y & & \\
 & h\swarrow & & \searrow j & \\
{\P}^3 & & & & ({\P}^3)^{\vee} \supset S \\
\end{array}$$
Moreover, each $\P^2$ fibre of $h:V\rightarrow\P^3$ is mapped to a hyperplane in $(\P^3)^{\vee}$ by $j$, and therefore each curve $Y_t$ is mapped (isomorphically) to a hyperplane section of $S\subset (\P^3)^{\vee}$ by $j$. Therefore the family of curves $Y/\P^3$ is a complete linear system of curves in a K3 surface, completing the proof.
\end{prf}

In Theorem~\ref{g=3nh} we assumed that every curve is non-hyperelliptic. Later we will investigate the hyperelliptic case, but here we show that the family cannot contain both non-hyperelliptic and hyperelliptic curves.

\begin{prp}
\label{notboth}
Let $Y\rightarrow\P^3$ be a flat family of integral Gorenstein curves of genus three whose compactified relative Jacobian $X=\overline{J}^d(Y/\P^3)$ is a Lagrangian fibration. If the family $Y/\P^3$ contains a non-hyperelliptic curve, then every curve in the family is non-hyperelliptic (and it follows that $X$ is a Beauville-Mukai system by Theorem~\ref{g=3nh}).
\end{prp}

\begin{prf}
Suppose that the family $Y/\P^3$ contains both non-hyperelliptic and hyperelliptic curves. Non-hyperelliptic curves are generic in the moduli space of genus three curves, so the hyperelliptic curves will occur in codimension at least one. The relative canonical map $Y\rightarrow V$ will map the non-hyperelliptic curves isomorphically to their images, but it will map the hyperelliptic curves two-to-one onto conics. Let $Y^{\prime}\subset V$ be the image of the canonical map, but with the conics replaced by `double conics' (i.e., given a non-reduced structure) so that $Y^{\prime}$ is a family of quartics over $\P^3$, given by the zero locus of a section of $\O(d+4,4)|_V$.
$$\begin{array}{ccccc}
Y & & \dashrightarrow & & Y^{\prime}\subset V \\
 & \pi\searrow & & \swarrow h & \\
 & & {\P}^3 & & \\
 \end{array}$$
By Lemma~\ref{preliminaries}~(7), $R^1\pi_*\O_Y\cong\Omega^1_{\P^3}$, and as before we can compute $R^1h_*\O_{Y^{\prime}}\cong\O_{\P^3}(-d-4)\otimes \Omega^1_{\P^3}$. Since $Y$ and $Y^{\prime}$ are birational, these sheaves must be isomorphic, and thus $d=-4$. As before, there is a quartic K3 surface $S\subset (\P^3)^{\vee}$ and $Y^{\prime}$ is a $\P^2$-bundle over $S$.

The curves in the family $Y^{\prime}/\P^3$ are the hyperplane sections of $S\subset (\P^3)^{\vee}$. In particular, $S$ contains hyperplane sections that are double conics. Without loss of generality, suppose that such a hyperplane is given by $z_0=0$ in homogeneous coordinates. Then the quartic defining $S$ looks like
$$F=Q(z_0,z_1,z_2,z_3)^2+z_0C(z_0,z_1,z_2,z_3)$$
where $Q$ is a quadric and $C$ is a cubic. Using subscripts to denote derivatives, we find that
\begin{eqnarray*}
F_0 & = & 2QQ_0+C+z_0C_0 \\
F_i & = & 2QQ_i+z_0C_i \qquad\qquad\mbox{for }i=1,2,3.
\end{eqnarray*}
In particular, $S$ has at least six singularities in the hyperplane $z_0=0$, where the quadric and cubic intersect. Consequently, $Y^{\prime}$ has $\P^2$ singular loci sitting above these singular points of $S$.

Now $Y$ is smooth by the last part of Lemma~\ref{preliminaries}. Since $Y$ and $Y^{\prime}$ are isomorphic on the locus of non-hyperelliptic curves, each $\P^2$ singular locus in $Y^{\prime}$ must be entirely contained in double conics. This means there is a family of double conic hyperplane sections of $S$ of dimension at least one, implying that $S$ is non-reduced at a generic point. This contradicts the fact that a generic hyperplane section of $S$ is a smooth non-hyperelliptic curve.
\end{prf}

\subsection{Genus four}

\begin{thm}
\label{g=4}
Let $Y\rightarrow\P^4$ be a flat family of integral Gorenstein non-hyperelliptic curves of genus four whose compactified relative Jacobian $X=\overline{J}^d(Y/\P^4)$ is a Lagrangian fibration. Then $X$ is a Beauville-Mukai integrable system, i.e., the family of curves is a complete linear system of curves in a K3 surface.
\end{thm}

\begin{prf}
By Theorem~\ref{rosenlicht}~(2), the relative canonical map
$$Y\rightarrow \P(\Omega^1_{\P^4})\cong V\subset\P^4\times (\P^4)^{\vee}$$
is a closed embedding, with each curve $Y_t$ embedded as a degree six curve in the corresponding $\P^3$ fibre of $\P(\Omega^1_{\P^4})$. We regard $Y$ as a subvariety of $V$.

Recall that a degree six curve $C$ in $\P^3$ is the complete intersection of a quadric and a cubic. Firstly, $C$ is contained in a unique quadric hypersurface: the short exact sequence
$$0\rightarrow\I_C(2)\rightarrow\O(2)\rightarrow\O(2)|_C\rightarrow 0$$
gives rise to a long exact sequence
$$0\rightarrow\H^0(\P^3,\I_C(2))\rightarrow\H^0(\P^3,\O(2))\rightarrow\H^0(C,\O(2)|_C)\rightarrow\H^1(\P^3,\I_C(2))\rightarrow\ldots$$
Now
$$\H^1(C,\O(2)|_C)\cong\H^0(C,\O(-2)|_C\otimes\omega_C)^{\vee}$$
vanishes so $h^0(C,\O(2)|_C)=9$ by Riemann-Roch, whereas $h^0(\P^3,\O(2))=10$. Therefore $h^0(\P^3,\I_C(2))\geq 1$ and $C$ is contained in at least one quadric. But if $C$ were contained in two distinct quadrics then it would have degree at most four. Note that this argument is valid whether $C$ is smooth or singular.

In the relative setting, we see that $Y$ will be contained in a unique hypersurface given by the zero locus of a section
$$q\in\H^0(\P(\Omega^1_{\P^4}),\O_{\P(\Omega^1_{\P^4})}(2)\otimes h^*\O_{\P^4}(d))=\H^0(V,\O(d+2,2)|_V)$$
for some integer $d$. By Lemma~\ref{aboutV}~(2), $d$ must be at least $-2$ for this space of sections to be non-trivial.

A similar computation as above shows that the space $\H^0(\P^3,\I_C(3))$ of cubics containing $C$ is five-dimensional. Four of these dimensions come from the following subspace
$$\H^0(\P^3,\I_C(2))\otimes\H^0(\P^3,\O(1))\hookrightarrow\H^0(\P^3,\I_C(3))$$
corresponding to cubics containing the quadric. Then $C$ is the complete intersection of the quadric and any cubic which does not contain the quadric. Because the cubic is only determined up to addition of a linear multiple of the quadric, extending this to the relative setting requires some care.

\vspace*{3mm}
\noindent
{\bf Claim:} There is a hypersurface given by the zero locus of a section
$$c\in\H^0(\P(\Omega^1_{\P^4}),\O_{\P(\Omega^1_{\P^4})}(3)\otimes h^*\O_{\P^4}(e))=\H^0(V,\O(e+3,3)|_V)$$
such that $Y$ is the complete intersection of the zero loci of $q$ and $c$.

\vspace*{3mm}
Because each degree six curve in $\P^3$ is the complete intersection of a quadric and a cubic, the required `relative cubic' exists locally in the base $\P^4$. In other words, there exists an open cover $\{U_i\}$ of $\P^4$ such that $Y|_{U_i}$ is the complete intersection of the quadric $\{q=0\}$ and a hypersurface in $V_i:=\P(\Omega^1_{\P^4})|_{U_i}$ given by the zero locus of a section
$$c_i\in\H^0(V_i,\O_{\P(\Omega^1_{\P^4})}(3)|_{V_i}).$$
On an intersection $V_i\cap V_j$, $c_i$ and $c_j$ must agree up to addition of a linear multiple of $q$ and up to an overall factor. The ``overall factor'' in this context is really the pull-back of a regular function $g_{ij}$ on $U_i\cap U_j\subset\P^4$, which is non-vanishing on $U_i\cap U_j$. Thus
$$c_i=h^*g_{ij}(c_j+l_{ij}q).$$
It is clear from the above equation that the $g_{ij}$ must satisfy $g_{ji}=g_{ij}^{-1}$ and $g_{ik}=g_{ij}g_{jk}$. In other words, they can be regarded as transition functions for a line bundle on $\P^4$, which can be absorbed into the definition of the $c_i$ by thinking of the relative cubics as sections
$$\tilde{c}_i\in\H^0(V_i,\O_{\P(\Omega^1_{\P^4})}(3)\otimes h^*\O_{\P^4}(e)|_{V_i}).$$
for some integer $e$. Then
$$\tilde{c}_i=\tilde{c}_j+\tilde{l}_{ij}q,$$
and to preserve degrees the `relative linear forms' $\tilde{l}_{ij}$ must be sections
$$\tilde{l}_{ij}\in\H^0(V_i\cap V_j,\O_{\P(\Omega^1_{\P^4})}(1)\otimes h^*\O_{\P^4}(e-d)|_{V_i\cap V_j}).$$
Now 
$$(\tilde{l}_{ij}+\tilde{l}_{jk}-\tilde{l}_{ik})q=(\tilde{c}_i-\tilde{c}_j)+(\tilde{c}_j-\tilde{c}_k)-(\tilde{c}_i-\tilde{c}_k)=0$$
and therefore $[\tilde{l}_{ij}]$ defines a class in
$$\H^1(\P(\Omega^1_{\P^4}),\O_{\P(\Omega^1_{\P^4})}(1)\otimes h^*\O_{\P^4}(e-d))=\H^1(V,\O(1+e-d,1)|_V).$$
By Lemma~\ref{aboutV}~(3), this cohomology group is trivial. Therefore we can write
$$\tilde{l}_{ij}=m_i-m_j$$
and replace $\tilde{c}_i$ by $\hat{c}_i:=\tilde{c}_i-m_iq$. On the overlap $V_i\cap V_j$,
$$\hat{c}_i-\hat{c}_j=\tilde{c}_i-m_iq-\tilde{c}_j+m_jq=\tilde{l}_{ij}q-m_iq+m_jq=0,$$
so the relative cubics $\hat{c}_i$ patch together to give a global section
$$c\in\H^0(\P(\Omega^1_{\P^4}),\O_{\P(\Omega^1_{\P^4})}(3)\otimes h^*\O_{\P^4}(e))=\H^0(V,\O(e+3,3)|_V),$$
proving the claim.

\vspace*{3mm}
By Lemma~\ref{aboutV}~(2), $e$ must be at least $-3$ for the above space of sections to be non-trivial.

\vspace*{3mm}
\noindent
{\bf Claim:} $d=-2$ and $e=-3$.

\vspace*{3mm}
By Lemma~\ref{preliminaries}~(7) we know that $R^1\pi_*\O_Y\cong\Omega^1_{\P^4}$. Thinking of $Y$ as the complete intersection of the hypersurfaces $\{q=0\}$ and $\{c=0\}$ in $V$, we compute this a second time using the projection $h:V\rightarrow\P^4$. We start with a resolution of $\O_Y$ by locally free sheaves:
$$0\rightarrow\O_{\P(\Omega^1_{\P^4})}(-5)\otimes h^*\O_{\P^4}(-d-e)\rightarrow
\begin{array}{c}
{\O}_{{\P}(\Omega^1_{{\P}^4})}(-2)\otimes h^*{\O}_{{\P}^4}(-d) \\
\oplus \\
{\O}_{{\P}(\Omega^1_{{\P}^4})}(-3)\otimes h^*{\O}_{{\P}^4}(-e)
\end{array}
\rightarrow\O_{\P(\Omega^1_{\P^4})}\rightarrow\O_Y\rightarrow 0$$
Since $h:V\rightarrow\P^4$ is a $\P^3$-bundle, and $\H^1$ and $\H^2$ vanish for all line bundles on $\P^3$, we deduce that $R^1h_*$ and $R^2h_*$ must vanish for all of the line bundles in the above resolution. This implies that the spectral sequence for $R^{\bullet}h_*$ degenerates and yields
\begin{eqnarray*}
R^1h_*{\O}_Y & \cong & R^3h_*\left({\O}_{{\P}(\Omega^1_{{\P}^4})}(-5)\otimes h^*{\O}_{{\P}^4}(-d-e)\right) \\
 & \cong & {\O}_{{\P}^4}(-d-e)\otimes R^3h_*\left({\O}_{{\P}(\Omega^1_{{\P}^4})}(-5)\right) \\
 & \cong & {\O}_{{\P}^4}(-d-e)\otimes \left(h_*\left({\O}_{{\P}(\Omega^1_{{\P}^4})}(5)\otimes\omega_h\right)\right)^{\vee} \\
 & \cong & {\O}_{{\P}^4}(-d-e)\otimes \left(h_*\left({\O}_{{\P}(\Omega^1_{{\P}^4})}(5)\otimes{\O}_{{\P}(\Omega^1_{{\P}^4})}(-4)\otimes h^*\omega_{{\P}^4}^{\vee}\right)\right)^{\vee} \\
 & \cong & {\O}_{{\P}^4}(-d-e)\otimes\omega_{{\P}^4}\otimes \left(h_*{\O}_{{\P}(\Omega^1_{{\P}^4})}(1)\right)^{\vee} \\
 & \cong & {\O}_{{\P}^4}(-d-e-5)\otimes \Omega^1_{{\P}^4}
 \end{eqnarray*}
where we have used relative Serre duality on the third line and the fact that
$$\omega_V\cong\omega_{\P^4\times(\P^4)^{\vee}}\otimes\O(1,1)|_V\cong\O(-4,-4)|_V\cong\O_{\P(\Omega^1_{\P^4})}(-4)$$
on the fourth line. Since $R^1h_*\O_Y\cong R^1\pi_*\O_Y\cong\Omega^1_{\P^4}$, we conclude that $d+e$ must equal $-5$. Since $d\geq -2$ and $e\geq -3$, we must have equality, proving the claim.

\vspace*{3mm}
We have proved that $Y\subset V$ is the complete intersection of the zero loci of sections $q$ and $c$ of $\O(0,2)|_V$ and $\O(0,3)|_V$, respectively. By Lemma~\ref{aboutV}~(2), these section lifts to sections of $\O(0,2)$ and $\O(0,3)$ on $\P^4 \times (\P^4)^{\vee}$, which define a degree six K3 surface $S$ in $(\P^4)^{\vee}$, the intersection of a quadric and a cubic. The projection $j:V\rightarrow (\P^4)^{\vee}$ expresses $Y$ as a $\P^3$-bundle over $S$. In particular, $S$ must be smooth since $Y$ is smooth by the last part of Lemma~\ref{preliminaries}. 
$$\begin{array}{ccccccc}
 & & V \supset Y & & \\
 & h\swarrow & & \searrow j & \\
{\P}^4 & & & & ({\P}^4)^{\vee} \supset S \\
\end{array}$$
Moreover, each $\P^3$ fibre of $h:V\rightarrow\P^4$ is mapped to a hyperplane in $(\P^4)^{\vee}$ by $j$, and therefore each curve $Y_t$ is mapped (isomorphically) to a hyperplane section of $S\subset (\P^4)^{\vee}$ by $j$. Therefore the family of curves $Y/\P^4$ is a complete linear system of curves in a K3 surface, completing the proof.
\end{prf}

\subsection{Genus five}

\begin{thm}
\label{g=5}
Let $Y\rightarrow\P^5$ be a flat family of integral Gorenstein non-hyperelliptic curves of genus five whose compactified relative Jacobian $X=\overline{J}^d(Y/\P^5)$ is a Lagrangian fibration. Then $X$ is a Beauville-Mukai integrable system, i.e., the family of curves is a complete linear system of curves in a K3 surface.
\end{thm}

\begin{prf}
By Theorem~\ref{rosenlicht}~(2), the relative canonical map
$$Y\rightarrow \P(\Omega^1_{\P^5})\cong V\subset\P^5\times (\P^5)^{\vee}$$
is a closed embedding, with each curve $Y_t$ embedded as a degree eight curve in the corresponding $\P^4$ fibre of $\P(\Omega^1_{\P^5})$. We regard $Y$ as a subvariety of $V$.

Recall that a degree eight curve $C$ in $\P^4$ is the complete intersection of a net of quadrics: the short exact sequence
$$0\rightarrow\I_C(2)\rightarrow\O(2)\rightarrow\O(2)|_C\rightarrow 0$$
gives rise to a long exact sequence
$$0\rightarrow\H^0(\P^4,\I_C(2))\rightarrow\H^0(\P^4,\O(2))\rightarrow\H^0(C,\O(2)|_C)\rightarrow\H^1(\P^4,\I_C(2))\rightarrow\ldots$$
Now
$$\H^1(C,\O(2)|_C)\cong\H^0(C,\O(-2)|_C\otimes\omega_C)^{\vee}$$
vanishes so $h^0(C,\O(2)|_C)=12$ by Riemann-Roch, whereas $h^0(\P^4,\O(2))=15$. Therefore $h^0(\P^4,\I_C(2))\geq 3$ and $C$ is contained in at least three independent quadrics. Since four independent quadrics would intersect in a zero-dimensional variety, $h^0(\P^4,\I_C(2))$ must be exactly $3$ and thus $C$ is contained in a net of quadrics. Note that this argument is valid whether $C$ is smooth or singular.


In the relative setting, we see that $Y$ will be the zero locus of a section
$$q\in\H^0(\P(\Omega^1_{\P^5}),\O_{\P(\Omega^1_{\P^5})}(2)\otimes h^*E)=\H^0(V,h^*E(2)\otimes j^*\O(2)|_V)$$
for some rank-three vector bundle $E$ on $\P^5$.

\vspace*{3mm}
\noindent
{\bf Claim:} The first Chern class of $E$ is $-6$.

\vspace*{3mm}
By Lemma~\ref{preliminaries}~(7) we know that $R^1\pi_*\O_Y\cong\Omega^1_{\P^5}$. Thinking of $Y$ as a subvariety of $V$, we compute this a second time using the projection $h:V\rightarrow\P^5$. We start with a resolution of $\O_Y$ by locally free sheaves:
$$0\rightarrow
\O_{\P(\Omega^1_{\P^5})}(-6)\otimes h^*\Lambda^3E^{\vee}\rightarrow
\O_{\P(\Omega^1_{\P^5})}(-4)\otimes h^*\Lambda^2E^{\vee}\hspace*{20mm}$$
$$\hspace*{20mm}\rightarrow\O_{\P(\Omega^1_{\P^5})}(-2)\otimes h^*E^{\vee}\rightarrow
\O_{\P(\Omega^1_{\P^5})}\rightarrow\O_Y\rightarrow 0$$
Since $h:V\rightarrow\P^5$ is a $\P^4$-bundle, and $\H^1$, $\H^2$, and $\H^3$ vanish for all line bundles on $\P^4$, we deduce that $R^1h_*$, $R^2h_*$, and $R^3h_*$ must vanish for all of the bundles in the above resolution (excluding $\O_Y$, of course). This implies that the spectral sequence for $R^{\bullet}h_*$ degenerates, and yields
\begin{eqnarray*}
R^1h_*{\O}_Y & \cong & R^4h_*\left({\O}_{{\P}(\Omega^1_{{\P}^5})}(-6)\otimes h^*\Lambda^3E^{\vee}\right) \\
 & \cong & \Lambda^3E^{\vee}\otimes R^4h_*\left({\O}_{{\P}(\Omega^1_{{\P}^5})}(-6)\right) \\
 & \cong & \Lambda^3E^{\vee}\otimes \left(h_*\left({\O}_{{\P}(\Omega^1_{{\P}^5})}(6)\otimes\omega_h\right)\right)^{\vee} \\
 & \cong & \Lambda^3E^{\vee}\otimes \left(h_*\left({\O}_{{\P}(\Omega^1_{{\P}^5})}(6)\otimes{\O}_{{\P}(\Omega^1_{{\P}^5})}(-5)\otimes h^*\omega_{{\P}^5}^{\vee}\right)\right)^{\vee} \\
 & \cong & \Lambda^3E^{\vee}\otimes\omega_{{\P}^5}\otimes \left(h_*{\O}_{{\P}(\Omega^1_{{\P}^5})}(1)\right)^{\vee} \\
 & \cong & \Lambda^3E^{\vee}(-6)\otimes \Omega^1_{{\P}^5}
 \end{eqnarray*}
where we have used relative Serre duality on the third line and the fact that
$$\omega_V\cong\omega_{\P^5\times(\P^5)^{\vee}}\otimes\O(1,1)|_V\cong\O(-5,-5)|_V\cong\O_{\P(\Omega^1_{\P^5})}(-5)$$
on the fourth line. Since $R^1h_*\O_Y\cong R^1\pi_*\O_Y\cong\Omega^1_{\P^5}$, we conclude that $\Lambda^3E^{\vee}(-6)$ is trivial. In other words, $\Lambda^3E\cong\O_{\P^5}(-6)$, proving the claim.

\vspace*{3mm}
\noindent
{\bf Case 1:} Suppose $E$ is a direct sum $\O(d_1)\oplus\O(d_2)\oplus\O(d_3)$ of line bundles on $\P^5$.

\vspace*{3mm}
If $E$ splits into a direct sum of line bundles, then the `relative net of quadrics' $q$ decomposes into three `relative quadrics'
$$q_i\in\H^0(\P(\Omega^1_{\P^5}),\O_{\P(\Omega^1_{\P^5})}(2)\otimes h^*\O_{\P^5}(d_i))=\H^0(V,\O(d_i+2,2)|_V)$$
for $i=1,2,3$. By Lemma~\ref{aboutV}~(2), each $d_i$ must be at least $-2$ for these spaces of sections to be non-trivial. But by the above claim, $c_1(E)=d_1+d_2+d_3=-6$, so we must have equality $d_1=d_2=d_3=-2$.

We have proved that $Y\subset V$ is the complete intersection of the zero loci of three sections $q_i$ of $\O(0,2)|_V$. By Lemma~\ref{aboutV}~(2), these section lifts to sections of $\O(0,2)$ on $\P^5 \times (\P^5)^{\vee}$, which define a degree eight K3 surface $S$ in $(\P^5)^{\vee}$, the intersection of a net of quadrics. The projection $j:V\rightarrow (\P^5)^{\vee}$ expresses $Y$ as a $\P^4$-bundle over $S$. In particular, $S$ must be smooth since $Y$ is smooth by the last part of Lemma~\ref{preliminaries}. 
$$\begin{array}{ccccccc}
 & & V \supset Y & & \\
 & h\swarrow & & \searrow j & \\
{\P}^5 & & & & ({\P}^5)^{\vee} \supset S \\
\end{array}$$
Moreover, each $\P^4$ fibre of $h:V\rightarrow\P^5$ is mapped to a hyperplane in $(\P^5)^{\vee}$ by $j$, and therefore each curve $Y_t$ is mapped (isomorphically) to a hyperplane section of $S\subset (\P^5)^{\vee}$ by $j$. Therefore the family of curves $Y/\P^5$ is a complete linear system of curves in a K3 surface, completing the proof in Case 1.

\vspace*{3mm}
\noindent
{\bf Case 2:} The general case.

\vspace*{3mm}
We will compute the degree of the discriminant locus $\Delta\subset\P^5$ parametrizing singular curves in the family $Y\rightarrow\P^5$ (equivalently, singular fibres of the Lagrangian fibration $X\rightarrow\P^5$), and then apply a theorem of the author~\cite{sawon08iv}. The degree of $\Delta$ can be computed using Chern classes. Recall that $Y$ is the zero locus of the section $q$, which can be regarded as a morphism
$$\O_{\P(\Omega^1_{\P^5})}(-2)\stackrel{q}{\longrightarrow} h^*E.$$
To find singularities in the fibres of $Y\rightarrow\P^5$, we consider
$$T_{\P(\Omega^1_{\P^5})/\P^5}(-2)\stackrel{dq}{\longrightarrow} h^*E$$
where on the left we have the relative tangent bundle of ${\P(\Omega^1_{\P^5})/\P^5}$ and the derivative $dq$ of $q$ is only in the fibre direction. Now $dq$ is a morphism from a rank-four bundle to a rank-three bundle, which will be generically surjective, so its kernel is a line bundle $L$ on $\P(\Omega^1_{\P^5})$. Moreover, $dq$ is locally given by a $4\times 3$ matrix, so it will drop rank on a  codimension two subset $\iota:\Sigma\hookrightarrow\P(\Omega^1_{\P^5})$. The cokernel of $dq$ will look like $\iota_*\mathcal{F}$, where $\mathcal{F}$ is a rank-one sheaf on $\Sigma$. We therefore have an exact sequence
$$0\rightarrow L\rightarrow T_{\P(\Omega^1_{\P^5})/\P^5}(-2)\stackrel{dq}{\longrightarrow} h^*E\rightarrow\iota_*\mathcal{F}\rightarrow 0.$$
We can use this to compute the first Chern class
$$c_1(L)=c_1(T_{\P(\Omega^1_{\P^5})/\P^5}(-2))-c_1(h^*E),$$
thereby determining $L$, and the second Chern class
$$c_2(\iota_*\mathcal{F})=[\Sigma].$$
Of course, the class $[Y]$ of $Y$ can be computed from the resolution
$$0\rightarrow
\O_{\P(\Omega^1_{\P^5})}(-6)\otimes h^*\Lambda^3E^{\vee}\rightarrow
\O_{\P(\Omega^1_{\P^5})}(-4)\otimes h^*\Lambda^2E^{\vee}\hspace*{20mm}$$
$$\hspace*{20mm}\rightarrow\O_{\P(\Omega^1_{\P^5})}(-2)\otimes h^*E^{\vee}\rightarrow
\O_{\P(\Omega^1_{\P^5})}\rightarrow\O_Y\rightarrow 0.$$
Combining these, we obtain the class $[Y\cap\Sigma]=[Y]\cap[\Sigma]$ of the locus of singular points of curves in the family $Y/\P^5$. If we restrict to a generic line $\ell\subset\P^5$, then the pencil of curves $Y|_{\ell}\rightarrow\ell\cong\P^1$ will contain exactly $\mathrm{deg}\Delta$ singular curves, each with a single simple node. The degree of $\Delta$ can therefore be computed from the bidegree of $[Y\cap\Sigma]$.

Now the only dependence on $E$ in the above calculation is through its Chern classes. Moreover, when we restrict to a line $\ell\subset\P^5$, only the rank and first Chern class of $E$ will appear, and these are the same as the rank and first Chern class of $\O(-2)\oplus\O(-2)\oplus\O(-2)$. But by Case 1, we know that we get a Beauville-Mukai integrable system if $E\cong\O(-2)\oplus\O(-2)\oplus\O(-2)$. Therefore, for general $E$ the degree of $\Delta$ will be the same as for the Beauville-Mukai system; this was computed to be $6n+18=48$ in Section~5 of~\cite{sawon08i}. Finally, the main result of~\cite{sawon08iv} states that (given a flat family $Y\rightarrow\P^n$ of integral Gorenstein curves whose compactified relative Jacobian $X\rightarrow\P^n$ is a Lagrangian fibration) $X$ is a Beauville-Mukai integrable system if $\mathrm{deg}\Delta>4n+20$. This is clearly satisfied in our case, since $n=5$ and  $\mathrm{deg}\Delta=48>4n+20=40$. This completes the proof.
\end{prf}

\begin{rmk}
We see, a postiori, that $E(2)$ must be trivial. It would be nice to have a direct proof of this fact. Of course, there exist indecomposable rank-three bundles on $\P^5$ with vanishing first Chern classes, as constructed by Horrocks~\cite{horrocks78}. But a calculation shows that $\H^0(V,h^*E(2)\otimes j^*\O(2)|_V)$ vanishes if $E(2)$ is a Horrocks bundle (at least for the ``parent bundle''), so they cannot be used to define families of genus five curves $Y/\P^5$.
\end{rmk}

\section{Hyperelliptic curves}

\subsection{Genus two}

The following result was originally proved by Markushevich~\cite{markushevich96}; our proof uses essentially the same ideas, with some simplifications.

\begin{thm}
\label{g=2}
Let $Y\rightarrow\P^2$ be a flat family of integral Gorenstein curves of genus two whose compactified relative Jacobian $X=\overline{J}^d(Y/\P^2)$ is a Lagrangian fibration. Then $X$ is a Beauville-Mukai integrable system, i.e., the family of curves is a complete linear system of curves in a K3 surface.
\end{thm}

\begin{prf}
By Theorem~\ref{rosenlicht}, integral Gorenstein curves of genus two must be hyperelliptic (note that there exists an integral non-Gorenstein non-hyperelliptic curve of genus two: it has a single singularity which looks locally like the intersection of the axes in $\C^3$, and its normalization is $\P^1$). Moreover, the relative canonical map
$$Y\rightarrow \P(\Omega^1_{\P^2})\cong V\subset\P^2\times (\P^2)^{\vee}$$
is two-to-one and onto, with each curve $Y_t$ a two-to-one cover of the corresponding $\P^1$ fibre of $\P(\Omega^1_{\P^2})$ branched over six points, counted with multiplicity. Thus the branch locus $B\subset V$ of the two-to-one map $Y\rightarrow V$ will be the zero locus of a section in
$$\H^0(\P(\Omega^1_{\P^2}),\O_{\P(\Omega^1_{\P^2})}(6)\otimes h^*\O_{\P^2}(2d))=\H^0(V,\O(2d+6,6)|_V)$$
for some integer $d$. By Lemma~\ref{aboutV}~(2), $d$ must be at least $-3$ for this space of sections to be non-trivial.

\vspace*{3mm}
\noindent
{\bf Claim:} $d=-3$.

\vspace*{3mm}
By Lemma~\ref{preliminaries}~(7) we know that $R^1\pi_*\O_Y=\Omega^1_{\P^2}$. On the other hand, $\pi:Y\rightarrow\P^2$ is the composition of the maps $g:Y\rightarrow\P(\Omega^1_{\P^2})$ and $h:\P(\Omega^1_{\P^2})\rightarrow\P^2$. The first direct image $R^1g_*$ must vanish, since $g$ is a two-to-one map, and
$$g_*\O_Y=\O_{\P(\Omega^1_{\P^2})}\oplus\O_{\P(\Omega^1_{\P^2})}\left(-\frac{1}{2}B\right)=\O_{\P(\Omega^1_{\P^2})}\oplus \left(\O_{\P(\Omega^1_{\P^2})}(-3)\otimes h^*\O_{\P^2}(-d)\right)$$
because the branch locus is $B$. To obtain $R^1\pi_*\O_Y$ we apply $R^1h_*$ to the above. The first term $R^1h_*\O_{\P(\Omega^1_{\P^2})}$ vanishes, and thus
\begin{eqnarray*}
R^1\pi_*{\O}_Y & \cong & R^1h_*\left(g_*{\O}_Y\right) \\
 & \cong & R^1h_*\left({\O}_{{\P}(\Omega^1_{{\P}^2})}(-3)\otimes h^*{\O}_{\P^2}(-d)\right) \\
 & \cong & {\O}_{{\P}^2}(-d)\otimes R^1h_*\left({\O}_{{\P}(\Omega^1_{{\P}^2})}(-3)\right) \\
 & \cong & {\O}_{{\P}^2}(-d)\otimes \left(h_*\left({\O}_{{\P}(\Omega^1_{{\P}^2})}(3)\otimes\omega_h\right)\right)^{\vee} \\
 & \cong & {\O}_{{\P}^2}(-d)\otimes \left(h_*\left({\O}_{{\P}(\Omega^1_{{\P}^2})}(3)\otimes{\O}_{{\P}(\Omega^1_{{\P}^2})}(-2)\otimes h^*\omega_{{\P}^2}^{\vee}\right)\right)^{\vee} \\
 & \cong & {\O}_{{\P}^2}(-d)\otimes\omega_{{\P}^2}\otimes \left(h_*{\O}_{{\P}(\Omega^1_{{\P}^2})}(1)\right)^{\vee} \\
 & \cong & {\O}_{{\P}^2}(-d-3)\otimes \Omega^1_{{\P}^2}
 \end{eqnarray*}
where we have used relative Serre duality on the fourth line and the fact that
$$\omega_V\cong\omega_{\P^2\times(\P^2)^{\vee}}\otimes\O(1,1)|_V\cong\O(-2,-2)|_V\cong\O_{\P(\Omega^1_{\P^2})}(-2)$$
on the fifth line. Since $R^1\pi_*\O_Y\cong\Omega^1_{\P^2}$, $d$ must equal $-3$, proving the claim.

\vspace*{3mm}
We have proved that $Y$ is a double cover of $V$ branched over the zero locus $B$ of a section of $\O(0,6)|_V$. By Lemma~\ref{aboutV}~(2), this section lifts to a section of $\O(0,6)$ on $\P^2 \times (\P^2)^{\vee}$, and we can define a K3 surface $S$ as the double cover of $(\P^2)^{\vee}$ branched over this sextic. The composition of maps
$$Y\stackrel{g}{\longrightarrow} V\stackrel{j}{\longrightarrow} (\P^2)^{\vee}$$
factors through $S$, and thereby expresses $Y$ as a $\P^1$-bundle over $S$. In particular, $S$ must be smooth since $Y$ is smooth by the last part of Lemma~\ref{preliminaries}. 
$$\begin{array}{ccccc}
 & & Y & & \\
 & & \downarrow g & \searrow & \\
 & & V & & S \\
 & h\swarrow & & \searrow j & \downarrow \\
{\P}^2 & & & & ({\P}^2)^{\vee} \\
\end{array}$$
Moreover, each $\P^1$ fibre of $h:V\rightarrow\P^2$ is mapped to a line in $(\P^2)^{\vee}$ by $j$, and the corresponding curve $Y_t$ is mapped (isomorphically) to the inverse image of that line in $S$. Therefore the family of curves $Y/\P^2$ is a complete linear system of curves in a K3 surface, completing the proof.
\end{prf}


\subsection{Reducible curves}

The canonical map of a genus three hyperelliptic curve $C$ will map $C$ two-to-one onto a conic in $\P^2$. Given a family of conics, we naturally expect there to be reducible conics (pairs of lines) in codimension one; a two-to-one cover of a pair of lines will likewise be reducible. Consequently, we expect that a family $Y/\P^3$ of hyperelliptic curves of genus three will contain reducible curves, and it is no longer reasonable to include integrality of the curves as a hypothesis.

A survey of different approaches to compactifying the Jacobian of a reducible curve can be found in Caporaso~\cite{caporaso08}. Many authors consider only curves with at worst ordinary double points, which is insufficient for our purposes. The most versatile approach seems to be to define the compactified Jacobian as a certain connected component of the Simpson moduli space of stable sheaves on the curve~\cite{simpson94}. When the curve has only ordinary double points, this compactification of the Jacobian is isomorphic to those constructed by Oda and Seshardi~\cite{os79} and Caporaso~\cite{caporaso94}, as was shown by Alexeev~\cite{alexeev04}. But it also applies to curves with more general singularities; for example, L{\' o}pez-Martin~\cite{lopez-martin05} used this approach to describe compactified Jacobians of tree-like curves and singular genus one curves of Kodaira type $I_m$, $III$, and $IV$. 

\begin{dfn}
Let $C$ be a reduced curve of arithmetic genus $n$, with polarization $H$. The Simpson Jacobian $J^d_sC$ is the moduli space of invertible sheaves on $C$ of degree $d$ which are stable with respect to the polarization $H$. The compactified (Simpson) Jacobian $\overline{J}^dC$ is the moduli space of stable rank-one (on every irreducible component) torsion-free sheaves on $C$ of degree $d$.
\end{dfn}

\begin{rmk}
The Hilbert polynomial $\chi(L\otimes\O_C(mH))$ of a rank-one torsion-free sheaf $L$ on $C$ is a linear function $am+b$. Simpson defined the slope of $L$ to be $\mu_H(L):=b/a$. Then $L$ is stable if its slope is greater than the slope of every proper non-trivial subsheaf. Similarly, $L$ is semi-stable if its slope is not less than the slope of every proper non-trivial subsheaf.
\end{rmk}

\begin{rmk}
By carefully choosing the degree $d$ and the polarization $H$, we can usually ensure that every semi-stable sheaf is automatically stable. If this were not the case, it would be natural to enlarge the above moduli spaces to include (equivalence classes of) semi-stable sheaves.
\end{rmk}

We illustrate this definition with some examples.

\begin{exm}
Let $C$ be a curve of arithmetic genus $n$ consisting of two smooth irreducible curves $C_1$ and $C_2$, touching at a pair of points. Thus the normalization $\tilde{C}$ is the disjoint union $C_1\sqcup C_2$,  there are points $p_1$ and $q_1$ in $C_1$ and $p_2$ and $q_2$ in $C_2$, and $C$ is given by joining $p_1$ to $p_2$ and $q_1$ to $q_2$, producing a pair of simple nodes. If $C_1$ has genus $n_1$ and $C_2$ has genus $n_2$, then $n=n_1+n_2+1$.

A degree $d$ line bundle $L$ on $C$ will consist of a degree $d_1$ line bundle $L_1$ on $C_1$ and a degree $d_2=d-d_1$ line bundle $L_2$ on $C_2$, plus {\em choices\/} of isomorphisms $(L_1)_{p_1}\cong (L_2)_{p_2}$ and $(L_1)_{q_1}\cong (L_2)_{q_2}$ up to an overall rescaling by $\C^*$. This yields an exact sequence
$$0\rightarrow\C^*\rightarrow J^dC\rightarrow\bigsqcup_{d_1+d_2=d}J^{d_1}C_1\times J^{d_2}C_2\rightarrow 0.$$
Clearly having infinitely many connected components is unmanageable, which is one reason to introduce the stability condition. Suppose that the polarization $H$ has degree $h_1$ on $C_1$ and degree $h_2$ on $C_2$. The Hilbert polynomial of $L$ is
$$\chi(L\otimes\O(mH))=(h_1+h_2)m+d+1-n$$
and its slope is
$$\mu_H(L)=\frac{d+1-n}{h_1+h_2}.$$
Now there is a surjection $L\rightarrow L\otimes\O_{C_2}=L_2$ whose kernel
$$L\otimes\mathcal{I}_{C_2\subset C}\cong L_1\otimes\mathcal{I}_{\{p_1,q_1\}\subset C_1}\cong L_1(-p_1-q_1)$$
has Hilbert polynomial
$$\chi(L_1(-p_1-q_1)\otimes\O(mH))=h_1m+d_1-2+1-n_1$$
and slope
$$\mu_H(L_1(-p_1-q_1))=\frac{d_1-1-n_1}{h_1}.$$
Since $L_1(-p_1-q_1)$ is a subsheaf of $L$, stability of $L$ implies
$$\frac{d_1-1-n_1}{h_1}<\frac{d+1-n}{h_1+h_2}.$$
Similarly, $L_2(-p_2-q_2)\subset L$ and stability of $L$ implies
$$\frac{d_2-1-n_2}{h_2}<\frac{d+1-n}{h_1+h_2}.$$
In this example, the above two inequalities are sufficient to establish the stability of $L$ (see Lemmas~3.3 and 3.4 of L{\'o}pez-Martin~\cite{lopez-martin05}).

For reasons which will become apparent shortly, let's assume that $n_1$ and $n_2$ are both positive and $H$ has degree $h_1=2n_1$ on $C_1$ and degree $h_2=2n_2$ on $C_2$. Since $d=d_1+d_2$ and $n=n_1+n_2+1$, the above inequalities become
$$\left|d_1-\left(\frac{n_1}{n_1+n_2}\right)d\right|<1\qquad\mbox{and}\qquad\left|d_2-\left(\frac{n_2}{n_1+n_2}\right)d\right|<1.$$

\noindent
{\bf Case 1:\/} If $n_1d/(n_1+n_2)$ is {\em not\/} an integer (equivalently, if $n_2d/(n_1+n_2)$ is {\em not\/} an integer) then there are two solutions which we will write as
$$(d_1,d_2)=(e_1+1,e_2)\qquad\mbox{and}\qquad (e_1,e_2+1).$$
The Simpson Jacobian of stable invertible sheaves on $C$ fits into the exact sequence
$$0\rightarrow\C^*\rightarrow J^d_sC\rightarrow J^{e_1+1}C_1\times J^{e_2}C_2\bigsqcup J^{e_1}C_1\times J^{e_2+1}C_2 \rightarrow 0.$$

\noindent
{\bf Case 2:\/} If $n_1d/(n_1+n_2)$ is an integer (equivalently, if $n_2d/(n_1+n_2)$ is an integer) then there is one solution
$$(d_1,d_2)=(e_1,e_2):=\left(\left(\frac{n_1}{n_1+n_2}\right)d,\left(\frac{n_2}{n_1+n_2}\right)d\right).$$
The Simpson Jacobian of stable invertible sheaves on $C$ fits into the exact sequence
$$0\rightarrow\C^*\rightarrow J^d_sC\rightarrow J^{e_1}C_1\times J^{e_2}C_2\rightarrow 0.$$

We compactify $J^d_sC$ by adding stable rank-one torsion-free sheaves. In the first case, $J^d_sC$ consists of two $\C^*$-bundles, over the $(n-1)$-dimensional abelian varieties $J^{e_1+1}C_1\times J^{e_2}C_2$ and $J^{e_1}C_1\times J^{e_2+1}C_2$, respectively. These are compactified by adding zero and infinity sections, $s^1_0$ and $s^1_{\infty}$, respectively, $s^2_0$ and $s^2_{\infty}$, to produce $\P^1$-bundles $P_1$ and $P_2$. The disjoint union $P_1\sqcup P_2$ is the normalization of the compactified Jacobian; it can be regarded as the scheme representing the Presentation Functor (see Altman and Kleiman~\cite{ak90}), or more precisely, an analogue of the Presentation Functor in the reducible case.

The compactified Jacobian $\overline{J}^dC$ itself is obtained by gluing $s^1_{\infty}$ to $s^2_0$ and  $s^2_{\infty}$ to $s^1_0$. Topologically $\overline{J}^dC$ looks like an $I_2$-bundle over an $(n-1)$-dimensional abelian variety. However, the zero and infinity sections are identified via the isomorphisms
\begin{eqnarray*}
s^1_{\infty}\cong J^{e_1+1}C_1\times J^{e_2}C_2 & \longrightarrow & J^{e_1}C_1\times J^{e_2+1}C_2\cong s^2_0 \\
(L_1,L_2) & \longmapsto & (L_1(-q_1),L_2(q_2))
\end{eqnarray*}
and
\begin{eqnarray*}
s^2_{\infty}\cong J^{e_1}C_1\times J^{e_2+1}C_2 & \longrightarrow & J^{e_1+1}C_1\times J^{e_2}C_2\cong s^1_0 \\
(L_1,L_2) & \longmapsto & (L_1(p_1),L_2(-p_2)).
\end{eqnarray*}
Consequently, the complex structure is not that of an $I_2$-bundle.

We will avoid the second case, which is more complicated. Although $J^d_sC$ now consists of a single $\C^*$-bundle over an $(n-1)$-dimensional abelian variety, it is not possible to compactify it without adding strictly semi-stable sheaves, such as invertible sheaves with $(d_1,d_2)=(e_1+1,e_2-1)$ or $(e_1-1,e_2+1)$.
\end{exm}

The compactified Jacobian in the next example was described by L{\'o}pez-Martin~\cite{lopez-martin05} in the case when the genera $n_1$ and $n_2$ are both zero (i.e., $C$ is a singular genus one curve of Kodaira type $III$). The higher genera cases can be deduced from this by generalizing the Presentation Functor~\cite{ak90} to reducible curves.

\begin{exm}
Let $C$ be a curve of arithmetic genus $n$ consisting of two smooth irreducible curves $C_1$ and $C_2$, touching tangentially at a point. Thus the normalization $\tilde{C}$ is the disjoint union $C_1\sqcup C_2$, there are points $p_1$ in $C_1$ and $p_2$ in $C_2$, and $C$ is given by joining $p_1$ to $p_2$ producing a tacnode. If $C_1$ has genus $n_1$ and $C_2$ has genus $n_2$, then $n=n_1+n_2+1$.

Note that $C$ is the limit of the previous example, when we let $q_1\rightarrow p_1$ and $q_2\rightarrow p_2$. In this sense, the relation between these two examples is similar to the relation between an irreducible curve with a node and an irreducible curve with a cusp, and one will notice the  similarities between their compactified Jacobians (as described in Subsection~2.1 of~\cite{sawon08iv}).

A degree $d$ line bundle $L$ on $C$ will consist of degree $d_1$ and $d_2$ line bundles on $C_1$ and $C_2$, respectively, plus a {\em choice\/} of isomorphism $(L_1)_{p_1}\cong (L_2)_{p_2}$ and compatible isomorphism of jet bundles $(\mathrm{Jet}^1L_1)_{p_1}\cong (\mathrm{Jet}^1L_2)_{p_2}$ up to an overall rescaling by $\C^*$. This yields an exact sequence
$$0\rightarrow\C\rightarrow J^dC\rightarrow\bigsqcup_{d_1+d_2=d}J^{d_1}C_1\times J^{d_2}C_2\rightarrow 0.$$
As in the previous example, one can determine the conditions on the degrees $d_1$ and $d_2$ in order for $L$ to be stable. If $n_1$ and $n_2$ are both positive and $H$ has degree $h_1=2n_1$ on $C_1$ and $h_2=2n_2$ on $C_2$, then $L$ is stable if and only if
$$\left|d_1-\left(\frac{n_1}{n_1+n_2}\right)d\right|<1\qquad\mbox{and}\qquad\left|d_2-\left(\frac{n_2}{n_1+n_2}\right)d\right|<1.$$
Once again, there are two cases. We will only consider the first case, when $n_1d/(n_1+n_2)$ and $n_2d/(n_1+n_2)$ are {\em not\/} integers, which yields two solutions
$$(d_1,d_2)=(e_1+1,e_2)\qquad\mbox{and}\qquad (e_1,e_2+1).$$
The Simpson Jacobian of stable invertible sheaves on $C$ fits into the exact sequence
$$0\rightarrow\C\rightarrow J^d_sC\rightarrow J^{e_1+1}C_1\times J^{e_2}C_2\bigsqcup J^{e_1}C_1\times J^{e_2+1}C_2 \rightarrow 0.$$

Thus $J^d_sC$ consists of two $\C$-bundles, over the $(n-1)$-dimensional abelian varieties $J^{e_1+1}C_1\times J^{e_2}C_2$ and $J^{e_1}C_1\times J^{e_2+1}C_2$, respectively, which we compactify by adding infinity sections $s^1_{\infty}$ and $s^2_{\infty}$, thereby producing $\P^1$-bundles $P_1$ and $P_2$. The disjoint union $P_1\sqcup P_2$ is the normalization of the compactified Jacobian; it can be regarded as the scheme representing the Presentation Functor corresponding to the partial normalization $C_1*C_2\rightarrow C$, where $C_1*C_2$ denotes the curve obtained by joining $p_1\in C_1$ to $p_2\in C_2$ producing a simple node (rather than a tacnode).

The compactified Jacobian $\overline{J}^dC$ itself consists of $P_1$ and $P_2$ touching `tangentially' along $s^1_{\infty}$ and $s^2_{\infty}$. In other words, we identify $s^1_{\infty}$ with $s^2_{\infty}$ using the isomorphism
\begin{eqnarray*}
s^1_{\infty}\cong J^{e_1+1}C_1\times J^{e_2}C_2 & \longrightarrow & J^{e_1}C_1\times J^{e_2+1}C_2\cong s^2_{\infty} \\
(L_1,L_2) & \longmapsto & (L_1(-p_1),L_2(p_2))
\end{eqnarray*}
{\em and\/} we identify a vector field $v^1$ along $s^1_{\infty}$ with a vector field $v^2$ along $s^2_{\infty}$. The vector field $v^1$ points out of $s^1_{\infty}$, and if we project down to $J^{e_1+1}C_1\times J^{e_2}C_2$, we get a vector field $d\pi(v^1)$ which at each point $(q_1,q_2)\in J^{e_1+1}C_1\times J^{e_2}C_2$ points in the direction of the curve $C_1$, if we take some Abel embedding
$$C_1\hookrightarrow J^{e_1+1}C_1\cong J^{e_1+1}C_1\times \{q_2\}\subset J^{e_1+1}C_1\times J^{e_2}C_2$$
and then translate it so that $p_1\in C_1$ is mapped to $(q_1,q_2)\in J^{e_1+1}C_1\times J^{e_2}C_2$. In particular, $v^1$ is {\em not\/} tangent to the $\P^1$ fibres of $P_1$. The vector field $v^2$ is defined in a similar way.

Gluing $s^1_{\infty}$ to $s^2_{\infty}$ and identifying $v^1$ with $v^2$ creates a family of tacnodes along the image $\mathrm{Sing}\overline{J}^dC$ of $s^1_{\infty}$ (and $s^2_{\infty}$): locally $\overline{J}^dC$ looks like the product of a tacnode and $\C^{n-1}$. Let $(z_1,z_2,...,z_n)$ be local coordinates on $P_1$, with $s^1_{\infty}$ given by $z_1=0$ and $v^1$ given by $\frac{\partial\phantom{z}}{\partial z_1}+\frac{\partial\phantom{z}}{\partial z_2}$. Similarly, let $(w_1,w_2,...,w_n)$ be local coordinates on $P_2$, with $s^2_{\infty}$ given by $w_1=0$ and $v^2$ given by $\frac{\partial\phantom{w}}{\partial w_1}+\frac{\partial\phantom{www}}{\partial w_{n_1+2}}$. Locally, the map to $\overline{J}^dC$ is given by
$$f(z_1,z_2,\ldots,z_n)=(z_1^2,z_1^4,z_2-z_1,z_3,\ldots,z_n)\in\{(x,y,z,\ldots )\in\C^{n+1}|y^2=x^4\}$$
on $P_1$ and by
$$f(w_1,w_2,\ldots,w_n)=(w_1^2,w_1^4,w_2,\ldots,w_{n_1+2}-w_1,\ldots,w_n)\in\{(x,y,z,\ldots )\in\C^{n+1}|y^2=x^4\}$$
on $P_2$. Note that pairs of $\P^1$ fibres in $P_1$ and $P_2$ combine to give curves $\P^1*\P^1$ in $\overline{J}^dC$ with simple nodes, rather than with tacnodes. Specifically, take a $\P^1$ fibre of $P_1$  given by keeping $z_2,\ldots,z_n$ constant and a $\P^1$ fibre of $P_2$ given by keeping $w_2,\ldots,w_n$ constant. If
$$(z_2,\ldots,z_n)=(w_2,\ldots,w_n)$$
then the image of these rational curves in $\overline{J}^dC$ will meet and produce a nodal curve $\P^1*\P^1$, rather than a singular genus one curve of Kodaira type $III$.
\end{exm}

Next we want to extend the preliminary results of Lemma~\ref{preliminaries} to the reducible case. We need some restrictions on the curves. The following definitions are from Catanese~\cite{catanese82}.

\begin{dfn}
A reduced curve $C$ is said to be canonically positive if it is Gorenstein and if for every irreducible component $C_i$ of $C$ the degree of $\omega_C|_{C_i}=\omega_C\otimes\O_{C_i}$ is positive.
\end{dfn}

\begin{rmk}
Note that $C$ is Gorenstein if and only if $\omega_C$ is an invertible sheaf, which is needed to define the degree. By Lemma~1.12~\cite{catanese82}, $\mathcal{I}_{C-C_i}\otimes\O_{C_i}$ is invertible,
$$\omega_C|_{C_i}=\omega_{C_i}\otimes (\mathcal{I}_{C-C_i}\otimes\O_{C_i})^{-1},$$
and
$$\mathrm{deg}\omega_C|_{C_i}=\mathrm{deg}\omega_{C_i}+C_i\cdot (C-C_i).$$
Therefore the degree of $\omega_C|_{C_i}$ is only non-positive if $C_i$ is a rational curve that intersects the rest of $C$ in at most two points (counted with multiplicity) or an isolated elliptic curve. Rational components intersecting the rest of $C$ in one or two points are `redundant' in two different senses: they are contracted by the canonical map, and collapsing them does not change the Jacobian or compactified Jacobian of the curve (though in a family, collapsing some components could create singularities in the total space). Therefore, we can say that a connected curve $C$ is canonically positive if it does not contain redundant rational components.
\end{rmk}

\begin{rmk}
In Caporaso, Coelho, and Esteves~\cite{cce08} canonically positive curves are called ``G-stable''. Note that a canonically positive curve $C$ with only simple nodes will be a stable curve.
\end{rmk}

\begin{rmk}
If $C$ is canonically positive, then $H=\omega_C$ is a natural polarization on $C$. Indeed, this is precisely the polarization we used in the two examples above.
\end{rmk}

\begin{dfn}
A curve $C$ is $2$-connected if for each decomposition $C=C_1\cup C_2$, with $C_1\cap C_2$ of dimension zero, we have $C_1\cdot C_2\geq 2$.
\end{dfn}

\begin{rmk}
If a connected curve $C$ is canonically positive, then
$$\mathrm{deg}\omega_C|_{C_i}=\mathrm{deg}\omega_{C_i}+C_i\cdot (C-C_i)\geq \mathrm{max}\{\mathrm{deg}\omega_{C_i},1\},$$
and only genus one components $C_i$ that intersect the rest of the curve in one point will be contracted by the canonical map of $C$. Such components cannot exist if $C$ is also $2$-connected. Therefore, if $C$ is both canonically positive and $2$-connected, then its canonical map does not contract any irreducible component of $C$.
\end{rmk}

\begin{exm}
One reason to avoid curves which are {\em not\/} $2$-connected is that their Simpson Jacobians could be empty. For example, suppose that $C$ consists of two smooth irreducible curves $C_1$ and $C_2$, of the same genus $n_1=n_2>0$, touching at a single point. The polarization $H=\omega_C$ has degree $h_1=2n_1-1$ on $C_1$ and degree $h_2=2n_2-1$ on $C_2$. In this case, the stability conditions for a degree $d=d_1+d_2$ line bundle on $C$ are
$$\left|d_1-\frac{1}{2}d\right|<\frac{1}{2}\qquad\mbox{and}\qquad\left|d_2-\frac{1}{2}d\right|<\frac{1}{2},$$
and there are no solutions if $d$ is odd.
\end{exm}

\begin{lem}
\label{Abel}
Let $C$ be a canonically positive $2$-connected curve. Then the Abel map, which takes a point $p\in C$ to $\mathfrak{m}^*_p\in\overline{J}^1C$, is an embedding.
\end{lem}

\begin{prf}
This lemma follows from results of Caporaso, Coelho, and Esteves. By Proposition~2 of~\cite{cce08}, if $C$ is ``free from separating nodes'' (i.e., $2$-connected) and ``G-stable'' (i.e., reduced, Gorenstein, and with ample dualizing sheaf; or in other words, canonically positive) then the Abel map embeds $C$ in the scheme parametrizing stable rank-one torsion-free sheaves on $C$ of degree one. Note that their definition of stability looks slightly different to ours, but in fact they are equivalent. We will just demonstrate a few of the simpler aspects of this result.



For a point $p$ in the smooth locus of $C$, let us verify that $\mathfrak{m}^*_p\cong\O_C(p)$ is stable. Let $C_i$ be an irreducible component of $C$ of arithmetic genus $n_i$ and with $C_i\cdot (C-C_i)=\delta_i$. Then the polarization $H=\omega_C$ has degree $h_i=2n_i-2+\delta_i>0$ on $C_i$. Summing over all irreducible components, we find
$$\sum h_i=\mathrm{deg}\omega_C=2n-2,$$
which can be used to show that the stability conditions for a degree $d=d_1+\ldots+d_k$ line bundle on $C$ are
$$d_i<\frac{1}{2}\delta_i+\left(\frac{h_i}{\sum h_i}\right)d.$$
Now $\O_C(p)$ has degree $d_i=1$ on the component $C_i$ containing $p$, and degree $d_j=0$ on all other components. These values satisfy the above inequalities, since $2$-connectivity of $C$ ensures that $\delta_j\geq 2$ for all $j$, and therefore $\O_C(p)$ is stable.

Next we show that if $p\neq q$ are points in the smooth locus of $C$, then $\O_C(p)\not \cong\O_C(q)$. If these line bundles are isomorphic, then $\O_C(p-q)$ has a nowhere vanishing section, i.e., there is a rational function $f$ on $C$ with a single pole at $p$ and a zero at $q$. If $p$ and $q$ lie in different components, then $f$ must vanish identically on the component containing $q$, contradicting the fact that $f$ came from an isomorphism $\O_C(p)\cong\O_C(q)$. If $p$ and $q$ lie in the same component $C_i$, then $f$ defines an isomorphism $C_i\cong\P^1$. Moreover, $f$ must be constant on all connected components of $C-C_i$. By $2$-connectivity, there is a connected component of $C-C_i$ which intersects $C_i$ in at least two points (counted with multiplicity); then $f$ takes the same value at these two points, contradicting the fact that $f:C_i\rightarrow\P^1$ is an isomorphism.

We've shown that the Abel map is birational onto its image in $\overline{J}^1C$. For the proof that it is actually an embedding we refer to~\cite{cce08}.
\end{prf}

\begin{exm}
Let $C$ consist of two smooth curves $C_1$ and $C_2$ touching at a pair of points, as in the first example of this subsection. The Abel map of $C$, restricted to $C_1-\{p_1,q_1\}$, gives a lift of the Abel embedding
$$C_1-\{p_1,q_1\}\subset C_1\hookrightarrow J^1C_1\cong J^1C_1\times\{\O_{C_2}\}\subset J^1C_1\times J^0C_2$$
of $C_1$ itself to the $\C^*$-bundle over $J^1C_1\times J^0C_2$ which makes up one connected component of $J^1_sC$. This lift extends to an embedding of $C_1$ into the $\P^1$-bundle $P_1$
$$\begin{array}{ccc}
 & & P_1 \\
 & \nearrow & \downarrow \\
C_1 & \hookrightarrow & J^1C_1\times J^0C_2.
\end{array}$$
Similarly, we have an embedding $C_2\hookrightarrow P_2$. When we glue $s^1_{\infty}\subset P_1$ to $s^2_0\subset P_2$, $q_1\in C_1$ is glued to $q_2\in C_2$. Similarly, when we glue $s^2_{\infty}\subset P_2$ to $s^1_0\subset P_1$, $p_2\in C_2$ is glued to $p_1\in C_1$. The result is an embedding of $C$ in $\overline{J}^1C$.
\end{exm}

\begin{lem}
\label{preliminaries2}
Let $Y\rightarrow B$ be a flat family of (reduced) canonically positive $2$-connected curves over a projective manifold $B$, such that the compactified relative Jacobian 
$$X=\overline{J}^d(Y/B)$$
is a Lagrangian fibration over $B$. (The relative compactified Jacobian $\overline{J}^d$ is defined in terms of the Simpson moduli space of stable sheaves~\cite{lopez-martin05, simpson94} with respect to the relative polarization $H=\omega_{Y/B}$. We implicitly assume that $d$ has been chosen so that every semi-stable sheaf is automatically stable, since otherwise $X$ would not be compact.) Then
\begin{enumerate}
\item$\hspace*{-2.5mm}^{\prime}$ every curve in the family $Y/B$ has arithmetic genus $n$,
\item$\hspace*{-2.5mm}^{\prime}$ the base $B$ is isomorphic to $\P^n$,
\item$\hspace*{-2.5mm}^{\prime}$ the generic curve in the family $Y/B$ is a smooth genus $n$ curve,
\item$\hspace*{-2.5mm}^{\prime}$ there is a hypersurface $\Delta\subset B$ parametrizing singular fibres of $\pi:X\rightarrow B$ (equivalently, singular curves in the family $Y/B$) and a curve above a generic point $\Delta$ will either be irreducible and contain a single simple node, or will consist of two smooth irreducible curves touching at a pair of points (as in the first example of this subsection),
\addtocounter{enumi}{1}
\item$\hspace*{-2.5mm}^{\prime}$ the first direct image sheaf $R^1\pi_*\O_X$ is isomorphic to $\Omega^1_{\P^n}$, where $\pi$ denotes the projection $X\rightarrow B\cong\P^n$,
\item$\hspace*{-2.5mm}^{\prime}$ if $d=1$, the first direct image sheaf $R^1\pi_*\O_Y$ is isomorphic to $\Omega^1_{\P^n}$, where $\pi$ denotes the projection $Y\rightarrow B\cong\P^n$.
\end{enumerate}
\end{lem}

\begin{prf}
The proofs of statements (1$^{\prime}$), (2$^{\prime}$), and (6$^{\prime}$) are the same as in the irreducible case; see Section~2.2 of Part~I~\cite{sawon08iv}. 

For statement (3$^{\prime}$), the generic fibre of the Lagrangian fibration $X/B$ is a smooth $n$-dimensional abelian variety. The only way for this to arise as the compactified Jacobian of a curve is if the curve $Y_t$ is tree-like. But since $Y_t$ is also $2$-connected, it must consist of a single smooth component.

The proof that $\Delta\subset B$ is a non-empty hypersurface is the same as in the irreducible case. Note that $X_t$ is smooth if and only if $Y_t$ is smooth, by the argument in the previous paragraph. Now suppose that $C=Y_t$ is a generic singular fibre of $Y/B$, with compactified Jacobian $\overline{J}^dC$. If the normalization $\tilde{C}$ of $C$ consists of a disjoint union
$$C_1\sqcup C_2\sqcup\ldots\sqcup C_k$$
of smooth curves of genus $n_1$, $n_2$, etc., then the same argument as in the irreducible case shows that the total genus $n_1+\ldots +n_k$ of $\tilde{C}$ must equal $n-1$. This means that the dual graph of $C$ has at most one loop. Using the facts that $C$ is both canonically positive and $2$-connected, we deduce that there are just four possibilities: either $C$ is irreducible with a single node or a single cusp, or $C$ consists of two smooth curves $C_1$ and $C_2$, with $n_1>0$ and $n_2>0$, touching at a pair of points or at a tacnode (as in the two examples at the start of this subsection). We already ruled out the cuspidal curve in the irreducible case, and the tacnodal curve can be eliminated by a similar argument. Namely, if such a curve occurred as a generic singular curve in codimension one, then the image of the pair of $\P^1$ fibres of $P_1$ and $P_2$ would be a characteristic $1$-cycle on $\overline{J}^dC$. However, we saw that these $\P^1$ fibres combine to give a curve $\P^1*\P^1$ in $\overline{J}^dC$, which is not an allowable characteristic $1$-cycle according to Hwang and Oguiso's classification~\cite{ho09}. This proves statement (4$^{\prime}$).

Finally, we consider statement (7$^{\prime}$). Since the curves in the family $Y/B$ all have arithmetic genus $n$, $R^1\pi_*\O_Y$ is locally free of rank $n$ with fibre
$$(R^1\pi_*\O_Y)_t \cong\H^1(Y_t,\O_{Y_t})$$
over the point $t\in B$. Since $d=1$, the relative Abel map gives a canonical map
$$Y\rightarrow X=\overline{J}^1(Y/B).$$
which is an embedding by Lemma~\ref{Abel}. We identify $Y$ with its image in $X$. The short exact sequence
$$0\rightarrow\mathcal{I}_Y\rightarrow\mathcal{O}_X\rightarrow\mathcal{O}_Y\rightarrow
0$$
yields the long exact sequence
$$R^1\pi_*\mathcal{I}_Y\rightarrow
R^1\pi_*\mathcal{O}_X\stackrel{\alpha}{\longrightarrow}R^1\pi_*\mathcal{O}_Y\rightarrow
R^2\pi_*\mathcal{I}_Y.$$
It remains to show that $\alpha$ is a surjection over $B\backslash\Delta_0$, where $\Delta_0\subset\Delta$ parametrizes non-generic singular fibres and is codimension two in $B$: since $R^1\pi_*\mathcal{O}_X$ and $R^1\pi_*\mathcal{O}_Y$ are both locally free of rank $n$ on $B$, surjectivity of $\alpha$ will imply that $\alpha$ is an isomorphism over $B\backslash\Delta_0$. Moreover, Hartogs' Theorem will then imply that $\alpha$ can be extended to a morphism
$$R^1\pi_*\mathcal{O}_X\longrightarrow R^1\pi_*\mathcal{O}_Y$$
over all of $B$, which must also be an isomorphism.


Let $t\in B\backslash\Delta_0$. If $t\not\in\Delta$, or if $t$ is a generic point of $\Delta$ and $Y_t$ is an irreducible curve with a single simple node, then
$$(R^1\pi_*\O_X)_t\cong\H^1(X_t,\O_{X_t})\cong\H^1(Y_t,\O_{Y_t})\cong(R^1\pi_*\O_Y)_t$$
by Lemma~7~\cite{sawon08iv}. It remains to show that the map
$$\H^1(X_t,\O_{X_t})\rightarrow\H^1(Y_t,\O_{Y_t})$$
induced from the Abel embedding $Y_t\hookrightarrow X_t$ is surjective when $C=Y_t$ consists of two smooth irreducible components $C_1$ and $C_2$ touching at a pair of points. The compactified Jacobian $X_t=\overline{J}^1C$ was described in the first example of this subsection. Let $g:C_1\sqcup C_2\rightarrow C$ and $g:P_1\sqcup P_2\rightarrow \overline{J}^1C$ be the normalizations of $C$ and $\overline{J}^1C$, respectively. We obtain the following commutative diagram
$$\begin{array}{ccccccccc}
0 & \rightarrow & \mathcal{O}_{\overline{J}^1C} & \rightarrow &
g_*\mathcal{O}_{P_1\sqcup P_2} & \rightarrow & \mathcal{G} & \rightarrow
& 0 \\
  & & \downarrow & & \downarrow & & \downarrow & & \\
0 & \rightarrow & \mathcal{O}_C & \rightarrow &
g_*\mathcal{O}_{C_1\sqcup C_2} & \rightarrow & \mathcal{G}^{\prime} & \rightarrow
& 0 \\
\end{array}$$
where $\mathcal{G}$ is supported on the singular locus $\mathrm{Sing}\overline{J}^1C$, which consists of two $(n-1)$-dimensional abelian varieties, and $\mathcal{G}^{\prime}$ is supported on the pair of nodes of $C$. The first vertical arrow comes from the Abel embedding of $C$ in $\overline{J}^1C$. The second vertical arrow comes from the embeddings $C_1\hookrightarrow P_1$ and $C_2\hookrightarrow P_2$, as described in the example following Lemma~\ref{Abel}. Finally, under the Abel embedding $C\hookrightarrow \overline{J}^1C$, the images of the two nodes of $C$ lie on $\mathrm{Sing}\overline{J}^1C$; the third vertical arrow comes from this inclusion. Taking cohomology we get
$$\begin{array}{ccccccccc}
0 & \rightarrow &
{\H}^0(\O_{\overline{J}^1C}) & \rightarrow &
{\H}^0(g_*\O_{P_1\sqcup P_2})\cong{\H}^0(\O_{P_1})\oplus\H^0(\O_{P_2}) & \rightarrow & {\H}^0(\mathcal{G}) & \rightarrow & \ldots \\
  & & \downarrow & & \downarrow & & \downarrow & & \\
0 & \rightarrow &
{\H}^0(\mathcal{O}_C) & \rightarrow &
{\H}^0(g_*\O_{C_1\sqcup C_2})\cong{\H}^0(\O_{C_1})\oplus\H^0(\O_{C_2})
& \rightarrow & {\H}^0(\mathcal{G}^{\prime}) & \rightarrow & \ldots \\
\end{array}$$
Both rows here look like
$$0\rightarrow\C\rightarrow\C^2\rightarrow\C^2\rightarrow\ldots$$
We can replace both $\H^0(\mathcal{G})$ and $\H^0(\mathcal{G}^{\prime})$ by the cokernels, isomorphic to $\C$, of the relevant maps and continue the long exact sequences
$$\begin{array}{ccccccccc}
{\C}_{(1)} & \stackrel{h_1}{\longrightarrow} &
{\H}^1(\O_{\overline{J}^1C}) & \stackrel{i_1}{\longrightarrow} &
{\H}^1(g_*\O_{P_1\sqcup P_2})\cong{\H}^1(\O_{P_1})\oplus{\H}^1(\O_{P_2})
& \stackrel{j_1}{\longrightarrow} & {\H}^1(\mathcal{G}) & \rightarrow & \ldots \\
\downarrow & & \downarrow & & \downarrow & & \downarrow & & \\
{\C}_{(2)} & \stackrel{h_2}{\longrightarrow} &
{\H}^1(\O_C) & \stackrel{i_2}{\longrightarrow} &
{\H}^1(g_*\O_{C_1\sqcup C_2})\cong{\H}^1(\O_{C_1})\oplus{\H}^1(\O_{C_2})
& \stackrel{j_2}{\longrightarrow} & {\H}^1(\mathcal{G}^{\prime})=0 & \rightarrow & \ldots \\
\end{array}$$
The first vertical arrow is an isomorphism and $h_1$ and $h_2$ are injections. Recall that $P_1$ and $P_2$ are $\P^1$-bundles over $J^1C_1\times J^0C_2$ and $J^0C_1\times J^1C_2$, respectively. The third vertical arrow comes from the compositions
$$\H^1(\O_{P_1})\cong\H^1(\O_{J^1C_1\times J^0C_2})\twoheadrightarrow\H^1(\O_{J^1C_1})\cong\H^1(\O_{C_1})$$
and
$$\H^1(\O_{P_2})\cong\H^1(\O_{J^0C_1\times J^1C_2})\twoheadrightarrow\H^1(\O_{J^1C_2})\cong\H^1(\O_{C_2}),$$
and is therefore surjective. Let $\beta\in\H^1(\O_C)$, and write
$$i_2(\beta)=(\beta_1,\beta_2)\in\H^1(\O_{C_1})\oplus\H^1(\O_{C_2}).$$
We can lift $(\beta_1,\beta_2)$ to
$$(\gamma_1,\gamma_2)\in\H^1(\O_{P_1})\oplus\H^1(\O_{P_2}).$$
We think of $\gamma_1$ and $\gamma_2$ as harmonic $(0,1)$-forms on $P_1$ and $P_2$, which therefore vanish in the $\P^1$-fibre directions and are translation invariant in the $(n-1)$-dimensional abelian variety directions. Now the map $j_1$ essentially takes $(\gamma_1,\gamma_2)$ to
$$(\gamma_1|_{s^1_{\infty}}-\gamma_2|_{s^2_0},\gamma_1|_{s^1_0}-\gamma_2|_{s^2_{\infty}}).$$
A priori, this won't vanish, but we can modify both $\gamma_1$ and $\gamma_2$ so that it does. To achieve this, we add an appropriate $\H^1(\O_{J^0C_2})$-component to
$$\gamma_1\in\H^1(\O_{P_1})\cong\H^1(\O_{J^1C_1\times J^0C_2})\cong\H^1(\O_{J^1C_1})\oplus\H^1(\O_{J^0C_2})$$
and an appropriate $\H^1(\O_{J^0C_1})$-component to
$$\gamma_2\in\H^1(\O_{P_2})\cong\H^1(\O_{J^0C_1\times J^1C_2})\cong\H^1(\O_{J^0C_1})\oplus\H^1(\O_{J^1C_2}).$$
In this way, we can create an element
$$(\gamma^{\prime}_1,\gamma^{\prime}_2)\in\H^1(\O_{P_1})\oplus\H^1(\O_{P_2}).$$
which still maps down to $(\beta_1,\beta_2)$, but is now in the kernel of $j_1$. Therefore we can lift $(\gamma^{\prime}_1,\gamma^{\prime}_2)$ to an element $\delta\in\H^1(\O_{\overline{J}^1C})$.

Let $\beta^{\prime}\in\H^1(\O_C)$ be the image of $\delta$ under the second vertical arrow. By construction $i_2(\beta-\beta^{\prime})=0$, so there exists $\epsilon\in\C_{(2)}$ such that $h_2(\epsilon)=\beta-\beta^{\prime}$. Under the isomorphism of the first vertical arrow, $\epsilon$ corresponds to $\kappa\in\C_{(1)}$. Finally, we let
$$\delta^{\prime}=\delta+h_1(\kappa)\in\H^1(\O_{\overline{J}^1C}).$$
Under the second vertical arrow, $\delta^{\prime}$ is mapped to
$$\beta^{\prime}+h_2(\epsilon)=\beta^{\prime}+\beta-\beta^{\prime}=\beta\in\H^1(\O_C).$$
This proves that the second vertical arrow
$$\begin{array}{c}
{\H}^1(\O_{\overline{J}^1C}) \\
\downarrow \\
{\H}^1(\O_C)
\end{array}$$
is surjective, completing the proof.
\end{prf}

\begin{rmk}
Let $S\rightarrow\P^1$ be an elliptic K3 surface containing a singular genus one curve $E_p$ of Kodaira type $III$ over the point $p\in\P^1$, i.e., $E_p$ consists of a pair of rational curves touching tangentially at a point. If $p_2,\ldots,p_n$ are generic points of $\P^1$, distinct from $p$ and each other, then the fibre of $\mathrm{Hilb}^nS\rightarrow\mathrm{Sym}^n\P^1$ over $\{p,p_2,\ldots,p_n\}$ is $E_p\times E_2\times\ldots \times E_n$. This is the kind of `tacnodal' fibre allowed by Hwang and Oguiso's classification.
\end{rmk}

\begin{rmk}
Lemma~\ref{preliminaries}~(5) has no analogue in the reducible case. Compactified Jacobians of reducible curves can change dramatically when we change the degree, so in general $\overline{J}^d(Y/B)$ and $\overline{J}^1(Y/B)$ will not be locally isomorphic as fibrations over $B$.
\end{rmk}

\begin{lem}
\label{smoothY}
Suppose that we are in the situation of Lemma~\ref{preliminaries2}, i.e., $Y\rightarrow B$ is a flat family of (reduced) canonically positive $2$-connected curves over a projective manifold $B$, such that the compactified relative Jacobian 
$$X=\overline{J}^d(Y/B)$$
is a Lagrangian fibration over $B$. Then the total space $Y$ of the family of curves is smooth.
\end{lem}

\begin{prf}
This is the analogue for reducible curves of the last statement of Lemma~\ref{preliminaries}. Our proof in the irreducible case used smoothness of the Abel map
$$\mathrm{Hilb}^m(Y/B)\rightarrow\overline{J}^m(Y/B)$$
for large degree $m$. For reducible curves, this map is not even well-defined: if all $m$ points lie on one component, the resulting rank-one torsion-free sheaf on the curve will not be stable. However, there is another approach which uses the bigraded Abel map of Esteves, Gagn{\'e}, and Kleiman~\cite{egk00}. They showed that for a flat family $Y/B$ of integral Gorenstein curves, the bigraded Abel map
$$Y\times_B J^{d-1}(Y/B)\rightarrow \overline{J}^d(Y/B)$$
is smooth. Since $\overline{J}^d(Y/B)$ is smooth by hypothesis, and $J^{d-1}(Y/B)$ is always smooth, it follows that $Y$ must be smooth.

In the reducible case, $\overline{J}^d(Y/B)$ becomes the moduli space of stable rank-one torsion-free sheaves. We {\em don't\/} add a stability condition to $J^{d-1}(Y/B)$, although this means that the bigraded Abel map will not be defined everywhere. Nevertheless, Esteves et al.'s proof shows that the map will be smooth wherever it is defined. Moreover, given $p\in Y$, there will always exist some family of degree $d-1$ line bundles on $Y/B$ such that the bigraded Abel map is defined in a neighbourhood of $p$. For example, suppose that $d=1$; then by Lemma~\ref{Abel} $\mathfrak{m}_q^*\otimes L$ is stable for any $q\in Y_t\subset Y$ and any line bundle $L$ on $Y_t$ which has degree zero on every component of $Y_t$. Therefore, as before, the smoothness of $J^{d-1}(Y/B)$ and $\overline{J}^d(Y/B)$ implies the smoothness of $Y$.
\end{prf}

\subsection{Genus three}

We adopt the following from Catanese~\cite{catanese82}, Definition~3.9 and Proposition~3.10.

\begin{dfn}
A reduced curve $C$ is hyperelliptic if there exist two points $p$ and $q$ in the smooth locus of $C$ such that $\H^0(\O_C(p+q))=2$.
\end{dfn}

\begin{prp}
Let $C$ be a reduced $2$-connected curve. Then $C$ is hyperelliptic if and only if the canonical map of $C$ is {\em not\/} birational.
\end{prp}

We can now prove the main result of this section.

\begin{thm}
\label{g=3h}
Let $Y\rightarrow\P^3$ be a flat family of (reduced) canonically positive $2$-connected hyperelliptic curves of genus three whose compactified Jacobian $X=\overline{J}^d(Y/\P^3)$ is a Lagrangian fibration, where $d$ is odd. Then $X$ is a Beauville-Mukai integrable system, i.e., the family of curves is a complete linear system of curves in a K3 surface.
\end{thm}

\begin{prf}
Tensoring with the relative polarization $H=\omega_{Y/\P^3}$, which has degree four on each curve, induces an isomorphism
$$\overline{J}^d(Y/\P^3)\cong\overline{J}^{d+4}(Y/\P^3).$$
We can also take dual sheaves fibre-wise, which induces an isomorphism
$$\overline{J}^d(Y/\P^3)\cong\overline{J}^{-d}(Y/\P^3).$$
So without loss of generality, we can assume $d=1$.

Recall that no irreducible component of a canonically positive $2$-connected curve is contracted by its canonical map. In fact, in genus three the canonical map either embeds the curve as a plane quartic (non-hyperelliptic case) or is a double cover of a plane conic (hyperelliptic case), by Lemma~2.7 of Catanese and Pignatelli~\cite{cp06}. Since all curves in the family $Y/\P^3$ are hyperelliptic, the relative canonical map
$$Y\rightarrow Z\subset\P(\Omega^1_{\P^3})\cong V\subset\P^3\times (\P^3)^{\vee}$$
is two-to-one onto its image $Z$, which is the zero locus of a section
$$c\in\H^0(\P(\Omega^1_{\P^3}),\O_{\P(\Omega^1_{\P^3})}(2)\otimes h^*\O_{\P^3}(d))=\H^0(V,\O(d+2,2)|_V)$$
for some integer $d$. By Lemma~\ref{aboutV}~(2), $d$ must be at least $-2$ for this space of sections to be non-trivial. Let $B\subset Z$ be the branch locus of the the two-to-one map $Y\rightarrow Z$.

Each curve $Y_t$ will be a two-to-one cover of its image, a conic, in the corresponding $\P^2$ fibre of $\P(\Omega^1_{\P^3})$ branched over eight points, counted with multiplicity. Let $D$ be a conic in $\P^2$. The short exact sequence
$$0\rightarrow\O(2)\rightarrow\O(4)\rightarrow\O(4)|_D\rightarrow 0$$
gives rise to the long exact sequence
$$0\rightarrow\H^0(\P^2,\O(2))\rightarrow\H^0(\P^2,\O(4))\rightarrow\H^0(D,\O(4)|_D)\rightarrow\H^1(\P^2,\O(2))=0.$$
We see that given any choice of eight points in the conic, there exists a quartic in $\P^2$ such that the eight points are given by the intersection of the quartic with the conic. This quartic is not unique: it is defined up to the addition of a quadric multiple of the quadric defining the conic $D$. We need to extend this to the relative setting.

\vspace*{3mm}
\noindent
{\bf Claim:\/} There is a hypersurface $B^{\prime}\subset V$ given by the zero locus of a section
$$q\in\H^0(\P(\Omega^1_{\P^3}),\O_{\P(\Omega^1_{\P^3})}(4)\otimes h^*\O_{\P^3}(2e))=\H^0(V,\O(2e+4,4)|_V)$$
such that $B=B^{\prime}\cap Z$.

\vspace*{3mm}
As explained above, the required quartic exists in each fibre, and indeed the required `relative quartic' exist locally in the base $\P^3$. The argument that we can patch these together into a global object is similar to the argument in the first claim of the proof of Theorem~\ref{g=4}. There exists an open cover $\{U_i\}$ of $\P^3$ such that $B|_{U_i}$ is the intersection of $Z|_{U_i}$ and a hypersurface in $V_i:=\P(\Omega^1_{\P^3})|_{U_i}$ given by the zero locus of a section
$$q_i\in\H^0(V_i,\O_{\P(\Omega^1_{\P^3})}(4)\otimes h^*\O_{\P^3}(2e)|_{V_i})$$
for some $e$, a priori a half integer. On an intersection $V_i\cap V_j$, $q_i$ and $q_j$ must agree up to addition of a quadric multiple of $c$ (note that we have already absorbed an overall factor into $q_i$ and $q_j$ by twisting by the line bundle $h^*\O_{\P^3}(2e)$). Therefore
$$q_i=q_j+l_{ij}c$$
where the relative quadrics $l_{ij}$ must be sections in
$$\H^0(V_i\cap V_j,\O_{\P(\Omega^1_{\P^3})}(2)\otimes h^*\O_{\P^3}(2e-d)|_{V_i\cap V_j}).$$
Now $[l_{ij}]$ defines a class in
$$\H^1(\P(\Omega^1_{\P^3}),\O_{\P(\Omega^1_{\P^3})}(2)\otimes h^*\O_{\P^3}(2e-d))=\H^1(V,\O(2+2e-d,2)|_V).$$
By Lemma~\ref{aboutV}~(3), this cohomology group is trivial. Writing $l_{ij}=m_i-m_j$, and replacing $q_i$ by $\tilde{q}_i:=q_i-m_ic$ yields relative quartics $\tilde{q}_i$ which agree on the overlaps $V_i\cap V_j$, and which therefore patch together to give a global section
$$q\in\H^0(\P(\Omega^1_{\P^3}),\O_{\P(\Omega^1_{\P^3})}(4)\otimes h^*\O_{\P^3}(2e))=\H^0(V,\O(2e+4,4)|_V),$$
proving the claim.

\vspace*{3mm}
By Lemma~\ref{aboutV}~(2), $e$ must be at least $-2$ for the above space of sections to be non-trivial.

\vspace*{3mm}
\noindent
{\bf Claim:\/} $d=-2$ and $e=-2$.

\vspace*{3mm}
By Lemma~\ref{preliminaries2}~(7$^{\prime}$) we know that $R^1\pi_*\O_Y=\Omega^1_{\P^3}$. On the other hand, $\pi:Y\rightarrow\P^3$ is the composition of the maps $g:Y\rightarrow Z\subset\P(\Omega^1_{\P^3})$ and $h:\P(\Omega^1_{\P^3})\rightarrow \P^3$. The first direct image $R^1g_*$ must vanish, since $g$ is a two-to-one map, and
\begin{eqnarray*}
g_*{\O}_Y & = & {\O}_Z\oplus{\O}_Z\left(-\frac{1}{2}B\right) \\
 & = & {\O}_Z\oplus{\O}_Z\left(-\frac{1}{2}B^{\prime}\cap Z\right) \\
 & = & {\O}_Z\oplus\left.\left({\O}_{\P(\Omega^1_{\P^3})}(-2)\otimes h^*{\O}_{\P^3}(-e)\right)\right|_Z
\end{eqnarray*}
because the branch locus is $B=B^{\prime}\cap Z$. We can resolve $\O_Z$ by locally free sheaves:
$$0\rightarrow\O_{\P(\Omega^1_{\P^3})}(-2)\otimes h^*\O_{\P^3}(-d)\rightarrow\O_{\P(\Omega^1_{\P^3})}\rightarrow\O_Z\rightarrow 0$$
This also yields a resolution of the second term in $g_*\O_Y$ by locally free sheaves:
$$0\rightarrow\O_{\P(\Omega^1_{\P^3})}(-4)\otimes h^*\O_{\P^3}(-d-e)\rightarrow\O_{\P(\Omega^1_{\P^3})}(-2)\otimes h^*\O_{\P^3}(-e)\hspace*{10mm}$$
$$\hspace*{40mm}\rightarrow\left.\left({\O}_{\P(\Omega^1_{\P^3})}(-2)\otimes h^*{\O}_{\P^3}(-e)\right)\right|_Z\rightarrow 0$$
Since $h:V\rightarrow\P^3$ is a $\P^2$-bundle, $R^1h_*$ and $R^2h_*$ must vanish for many of the terms in the above resolutions. We find that the long exact sequence obtained by applying $R^{\bullet}h_*$ to the resulting resolution of $g_*\O_Y$ by locally free sheaves yields
\begin{eqnarray*}
R^1h_*\left(g_*\O_Y\right) & \cong & R^2h_*\left({\O}_{{\P}(\Omega^1_{{\P}^3})}(-4)\otimes h^*{\O}_{{\P}^3}(-d-e)\right) \\
 & \cong & {\O}_{{\P}^3}(-d-e)\otimes R^2h_*\left({\O}_{{\P}(\Omega^1_{{\P}^3})}(-4)\right) \\
 & \cong & {\O}_{{\P}^3}(-d-e)\otimes \left(h_*\left({\O}_{{\P}(\Omega^1_{{\P}^3})}(4)\otimes\omega_h\right)\right)^{\vee} \\
 & \cong & {\O}_{{\P}^3}(-d-e)\otimes \left(h_*\left({\O}_{{\P}(\Omega^1_{{\P}^3})}(4)\otimes{\O}_{{\P}(\Omega^1_{{\P}^3})}(-3)\otimes h^*\omega_{{\P}^3}^{\vee}\right)\right)^{\vee} \\
 & \cong & {\O}_{{\P}^3}(-d-e)\otimes \omega_{{\P}^3}\otimes\left(h_*{\O}_{{\P}(\Omega^1_{{\P}^3})}(1)\right)^{\vee} \\
 & \cong & {\O}_{{\P}^3}(-d-e-4)\otimes\Omega^1_{{\P}^3}
\end{eqnarray*}
where we have used relative Serre duality on the third line and the fact that
$$\omega_V\cong\omega_{\P^3\times (\P^3)^{\vee}}\otimes\O(1,1)|_V\cong\O(-3,-3)|_V\cong\O_{\P(\Omega^1_{\P^3})}(-3)$$
on the fourth line. Since $R^1h_*\left(g_*\O_Y\right)\cong R^1\pi_*\O_Y\cong\Omega^1_{\P^3}$, we conclude that $d+e$ must equal $-4$. Since $d\geq -2$ and $e\geq -2$, we must have equality, proving the claim.

\vspace*{3mm}
We have proved that $Y$ is a double cover of the zero locus $Z\subset V$ of a section of $\O(0,2)|_V$ branched over the intersection $B$ of $Z$ with the zero locus $B^{\prime}$ of a section of $\O(0,4)|_V$. By Lemma~\ref{aboutV}~(2), these sections lift to sections of $\O(0,2)$ and $\O(0,4)$, respectively, on $\P^3\times (\P^3)^{\vee}$, and we can define a K3 surface $S$ as the double cover of the resulting quadric $W$ in $(\P^3)^{\vee}$ branched over the intersection of $W$ with the quartic. The composition of maps
$$Y\stackrel{g}{\longrightarrow} Z\stackrel{j}{\longrightarrow} W$$
factors through $S$, and thereby expresses $Y$ as a $\P^1$-bundle over $S$. In particular, $S$ must be smooth since $Y$ is smooth by Lemma~\ref{smoothY}. 
$$\begin{array}{ccccc}
 & & Y & & \\
 & & \downarrow g & \searrow & \\
 & & Z\subset V  & & S \\
 & h\swarrow & & \searrow j & \downarrow \\
{\P}^3 & & & & W\subset ({\P}^3)^{\vee}
\end{array}$$
Moreover, each $\P^2$ fibre of $h:V\rightarrow\P^3$ is mapped to a hyperplane in $(\P^3)^{\vee}$ by $j$, intersecting $W$ in a conic (the image of $Z_t$ under $j$), and the corresponding curve $Y_t$ is mapped (isomorphically) to the inverse image of this conic in $S$. Therefore the family of curves $Y/\P^3$ is a complete linear system of curves in a K3 surface, completing the proof.
\end{prf}

\begin{rmk}
Suppose that $d$ is even. If $Y/\P^3$ contains a curve $C$ consisting of two smooth curves $C_1$ and $C_2$ of genus one touching at a pair of points, then there will exist strictly semi-stable sheaves on $C$. In particular, $\overline{J}^d(Y/\P^3)$ will be non-compact. If $Y/\P^3$ is a complete linear system of hyperelliptic curves on a K3 surface $S$, then $\overline{J}^d(Y/\P^3)$ can be defined as a Mukai moduli space of stable sheaves on $S$. In this case, it is possible to vary the polarization of $S$ so that all semi-stable sheaves must be stable: the induced relative polarization $H$ of the family of curves $Y/\P^3$ will be such that $H$ has different degrees $h_1\neq h_2$ on $C_1$ and $C_2$, so that strictly semi-stable sheaves on $C$ no longer exist. For a generic polarization on $S$, $\overline{J}^d(Y/\P^3)$ will therefore be compact.

On the other hand, for a general family of (canonically positive $2$-connected hyperelliptic) curves $Y/\P^3$, the relative canonical sheaf $\omega_{Y/\P^3}$ is the only relative polarization guaranteed to exist, and we cannot necessarily vary it to eliminate strictly semi-stable sheaves.
\end{rmk}

\begin{flushleft}
Department of Mathematics\hfill sawon@math.unc.edu\\
University of North Carolina\hfill www.unc.edu/$\sim$sawon\\
Chapel Hill NC 27599-3250\\
USA\\
\end{flushleft}


\begin{thebibliography}{XXX}


\bibitem{alexeev04} V. Alexeev,
{\em Compactified Jacobians and Torelli map\/},
Publ. Res. Inst. Math. Sci. {\bf 40} (2004), no. 4, 1241--1265.



\bibitem{ak90} A. Altman and S. Kleiman,
{\em The presentation functor and the compactified Jacobian\/},
The Grothendieck Festschrift, Vol. I, 15--32, Progr. Math. {\bf 86},
Birkh{\"a}user, 1990.



\bibitem{beauville99} A. Beauville,
{\em Counting rational curves on K3 surfaces\/},
Duke Math. J. {\bf 97} (1999), no. 1, 99--108.







\bibitem{caporaso94} L. Caporaso,
{\em A compactification of the universal Picard variety over the moduli space of stable curves\/},
J. Amer. Math. Soc. {\bf 7} (1994), no. 3, 589–-660.

\bibitem{caporaso08} L. Caporaso,
{\em Compactified Jacobians, Abel maps, and theta divisors\/},
in Curves and abelian varieties, pp. 1--23, Contemp. Math. {\bf 465}, Amer. Math. Soc., Providence, RI, 2008.


\bibitem{cce08} L. Caporaso, J. Coelho, and E. Esteves,
{\em Abel maps of Gorenstein curves\/},
Rend. Circ. Mat. Palermo (2) {\bf 57} (2008), no. 1, 33--59.

\bibitem{catanese82} F. Catanese,
{\em Pluricanonical-Gorenstein-curves\/},
in Enumerative geometry and classical algebraic geometry (Nice, 1981), pp. 51--95, Progr. Math. {\bf 24}, Birkh{\"a}user Boston, 1982.

\bibitem{cp06} F. Catanese and R. Pignatelli,
{\em Fibrations of low genus I.\/},
Ann. Sci. {\' E}cole Norm. Sup. (4) {\bf 39} (2006), no. 6, 1011--1049.








\bibitem{egk00} E. Esteves, M. Gagn{\'e}, and S. Kleiman,
{\em Abel maps and presentation schemes\/},
Comm. Algebra {\bf 28} (2000), no. 12, 5961--5992.






\bibitem{gh78} P. Griffiths and J. Harris,
{\em Principles of algebraic geometry\/},
Pure and Applied Mathematics, Wiley, New York, 1978.











\bibitem{horrocks78} G. Horrocks,
{\em Examples of rank three vector bundles on five-dimensional projective space\/},
J. London Math. Soc. (2), {\bf 18} (1978), 15--27.








\bibitem{ho09} J.--M. Hwang and K. Oguiso,
{\em Characteristic foliation on the discriminantal hypersurface of a
  holomorphic Lagrangian fibration\/},
Amer. J. Math. {\bf 131} (2009), no. 4, 981--1007.





\bibitem{km09} S. Kleiman and R. Martins,
{\em The canonical model of a singular curve\/},
Geom. Dedicata {\bf 139} (2009), 139--166.






\bibitem{lopez-martin05} A.C. L{\'o}pez-Martin,
{\em Simpson Jacobians of reducible curves\/},
J. Reine Angew. Math. {\bf 582} (2005), 1--39.



\bibitem{markushevich96} D. Markushevich,
{\em Lagrangian families of Jacobians of genus 2 curves\/},
J. Math. Sci. {\bf 82} (1996), no. 1, 3268--3284.



\bibitem{matsushita99} D. Matsushita,
{\em On fibre space structures of a projective irreducible symplectic
manifold\/},
Topology {\bf 38} (1999), no. 1, 79--83.
Addendum, Topology {\bf 40} (2001), no. 2, 431--432.

\bibitem{matsushita00i} D. Matsushita,
{\em Equidimensionality of Lagrangian fibrations on holomorphic
symplectic manifolds\/},
Math. Res. Lett. {\bf 7} (2000), no. 4, 389--391.







\bibitem{mukai84} S. Mukai,
{\em Symplectic structure of the moduli space of simple sheaves on an
abelian or K3 surface\/},
Invent. Math. {\bf 77} (1984), 101--116.




\bibitem{os79} T. Oda and C.S. Seshadri,
{\em Compactifications of the generalized Jacobian variety\/},
Trans. Amer. Math. Soc. {\bf 253} (1979), 1--90.








\bibitem{rosenlicht52} M. Rosenlicht,
{\em Equivalence relations on algebraic curves\/},
Ann. Math. {\bf 56} (1952), 169--191.




\bibitem{sawon08i} J. Sawon,
{\em On the discriminant locus of a Lagrangian fibration\/},
Math. Ann. {\bf 341} (2008), no. 1, 201--221.


\bibitem{sawon08iv} J. Sawon,
{\em On Lagrangian fibrations by Jacobians I\/},
preprint {\bf arXiv:0803.1186v3}.



\bibitem{simpson94} C. Simpson,
{\em Moduli of representations of the fundamental group of a smooth projective variety I.\/},
Inst. Hautes {\' E}tudes Sci. Publ. Math. {\bf 79} (1994), 47--129.












\end{thebibliography}
\end{document}